# Delivery by Drones with Arbitrary Energy Consumption Models: A New Formulation Approach

Amir Ahmadi-Javid[1] and Mahla Meskar

*Department of Industrial Engineering & Management Systems, Amirkabir University of Technology, Tehran, Iran*

2022-11-01

**Abstract.** This paper presents a new approach for formulating the delivery problem by drones with *general* energy consumption models where the drones visit a set of places to deliver parcels to customers. Drones can perform multiple trips that start and end at a central depot while visiting several customers along their paths. The problem determines the routing and scheduling decisions of the drones in order to minimize the total transportation cost of serving customers. For the first time, the new formulation approach enables us to use the best available energy consumption model without the need of any extra approximations. Though the approach works in a very general setting including non-convex energy consumption models, it is also computationally efficient as the resulting optimization model has a linear relaxation. A numerical study on 255 benchmark instances with up to 50 customers and a specific energy function indicate that all the instances can be solved 20 times faster on average using the new formulation when compared to the best existing branch-and-cut algorithm. All the 15 benchmark instances with 50 customers are solved exactly, whereas none of them has been solved optimally before. Moreover, new instances with up to 150 customers are solved with small error bounds within a few hours. The new approach can be simply applied to consider the extra energy required when a drone needs to continue hovering until opening the delivery time window. It can also be applied to the case where the flight time is dependent on the drone's payload weight. Owing to the flexibility of the new approach, these challenging extensions are formulated as linear optimization models for the first time.

**Keywords:** Last-mile delivery; Unmanned Aerial Vehicle (UAV), Drone routing and scheduling; Energy consumption model; Vehicle routing problem (VRP)

---

[1] Corresponding author's email address: ahmadi_javid@aut.ac.ir

# 1. Introduction

Since COVID-19 epidemic spreads around the world, producers and traders have faced challenges in continuing their business by adhering to social distance. E-commerce, which has already experienced an exceptional growth, is now more commonly used. By rising e-commerce and delivery services in addition to the improvements in drones' technologies and the need for social distancing measures, using drones in last-mile delivery has grown significantly. The number of recent publications studying routing and scheduling of drones used for parcel delivery is rising fast.

Drones have special characteristics that make them different from typical ground vehicles. Unlike trucks, drones are not limited to roads and streets. They do not delay in traffic and can fly nonstop to their destinations. Therefore, it is considered that delivery of small packages via drones can potentially be faster and greener, as well as less expensive than ground vehicles (Carlsson & Song, 2017; Goodchild & Toy, 2018). On the other hand, they have significant restrictions on their maximum flight duration and payload capacity. They need to reload delivery products and recharge battery very fast. The battery consumption rate of a drone between two particular points depends on the amount of drone's payload weight. Additionally, other factors such as the wind direction and speed can affect the maximum flight range of a drone.

The greatest challenge of using drones in parcel delivery is how to ensure on-time and unfailing parcel delivery. Schedules of drone delivery services may fail as a result of a strong wind, rain, and underestimating the energy consumption of these vehicles. In fact, how a drone's battery is consumed, determines the key performance metrics of its fly range, cost, and emissions. To explore drones' capability to provide fast delivery, with low cost and emissions in delivery operations, one needs to formulate their energy consumption models precisely. Key factors affecting a drone's energy consumption are the drone's design (e.g., drone's weight and size, battery weight, lift-to-drag ratio, power transfer efficiency, size and number of rotors, avionics), environmental factors (e.g., gravity value, wind direction, and weather condition, air density and temperature), and controllable factors (e.g., value and angle of drone's speed and acceleration, and weight and size of the payload). Moreover, different drone motions such as hovering, steady flight, departure, and landing follow different energy consumption models (Zhang et al., 2020). Existing optimization models that consider a few number of key factors influencing drones' energy consumption become very complex, and consequently computationally intractable.

The main contribution of this paper is to present an efficiently solvable formulation of a drone routing problem. The new formulation approach can use the best available energy consumption model without the need of any further approximations. In the drone routing problem, a fleet of drones visit a set of locations to deliver parcels to the customers within predefined time windows. Drones can perform multiple trips where each trip starts and ends at the depot while visiting several customers along the way. The problem determines the routing and scheduling decisions of the drones in order to minimize the total transportation cost of serving customers. The new formulation can be solved by general optimization solvers. It needs much less computational times (20 times faster) for optimally solving the benchmark instances with linear and convex energy functions, when compared to a recent algorithm developed for solving them.

The remainder of the paper is organized as follows. Section 2 reviews the related literature, and states the main contributions of this study. Section 3 includes the formal statement and old formulation of the problem. Section 4 presents the new approach for formulating the problem. Section 5 evaluates the performance of the new solution method using an extensive numerical study. Finally, Section 6 provides the conclusions and future research opportunities.

## 2. Literature review

The number of publications studying drone routing is growing rapidly in recent years and various drone routing problems have been introduced and studied. In general, the types of the problems investigated in the literature can be classified according to the types of available vehicles and how they work together, the number of each vehicle type, the number of parcels each drone can carry, and the energy consumption model. In terms of the types of available vehicles and how they work together, the literature can be classified into three classes: (1) Drone-only fleet, (2) Truck-drone fleet where they work independently, and (3) a fleet of truck-drone tandems working in cooperation where drones can be carried by trucks; in the third category trucks work as range extenders for drones. The review papers by Macrina et al. (2020), Otto et al. (2018), Li et al. (2021), and Rojas Viloria et al. (2021) classified the literature based on the available vehicle types, the number of each vehicle type, the cooperation method, and the solution methods. Furthermore, Zhang et al. (2020) reviewed the key published energy consumption models for drones and discussed similarities and differences among them. In the following, a number of papers not covered in these surveys are briefly reviewed.

Roberti and Ruthmair (2021) proposed a new compact mixed-integer linear program (MILP) formulation for several TSP-D variants in which a truck and a drone cooperate to serve a set of customers. The formulation is based on timely synchronizing truck and drone flows. The authors designed a branch-and-price algorithm and solved instances with up to 39 customers to optimality. Kang and Lee (2022) presented an exact algorithm based on the Benders decomposition method for a drone-truck routing problem with a heterogeneous fleet. Mbiadou Saleu et al. (2022) studied a PDSTSP in which a truck and several drones work independently to deliver parcels to customers with the aim of minimizing the completion time; the problem extends the classic PDSTSP by considering several vehicles. The authors proposed an MILP model, a branch-and-cut algorithm, and a metaheuristic to solve the instances. Nguyen et al. (2022) studied a similar problem with the objective of minimizing the total transportation costs, and designed a heuristic Ruin-and-Recreate algorithm.

Baloch and Gzara (2020) designed a distribution network for an e-retailer under service-based competition and explored the advantages of using drones for parcel delivery considering technological limitations, regulations, and customer behavior. In another study, Chen et al. (2022) used a deep Q-learning method to learn the value of assigning a new customer to a drone (or a vehicle) or reject the delivery service in the context of a same-day delivery problem with vehicles and drones. Dayarian et al. (2020) suggested a new delivery system in which trucks are regularly resupplied by drones and measured the potential benefits of the newly introduced delivery system. Xia et al. (2021) suggested extending drones' battery lifespans through a blockchain-enabled drone sharing approach. In their problem, multiple operators are available and every operator owns a depot. The collaboration happens when a drone flies to other operators' depots to swap its battery. They considered no energy consumption function for drones and only set a maximum fly range constraint.

Considering the problem addressed in this paper, the rest of this section focuses on routing problems with a drone-only fleet in which an energy consumption model is considered and the aim is parcel delivery.

## 2.1. DRP literature

This subsection reviews the literature of delivery problems using a drone-only fleet with an energy consumption model, which are simply called Drone Routing Problems (DRPs) in the sequel.

Most studies on DRPs have tried to simplify evaluating drones' energy consumption and formulated a drone's power depletion explicitly as a linear or convex function of the drone's weight. Dorling et al. (2017) studied a DRP in which the energy consumption model was a linear function of the drone's battery and payload weight where the battery weight was the decision variable. They explored two variants of the problem; the first minimizes the total cost with a delivery time restriction, and the second minimizes the delivery completion time subject to a budget limit (no time window is imposed as controlling/minimizing the maximum completion time is considered instead). The authors presented an MILP formulation, and developed a simulated annealing heuristic to solve the problem instances. They tested the presented heuristic on small instances with 6–8 delivery points, and large scenarios with 125 or 500 delivery locations. Rabta et al. (2018) examined the problem of delivering parcels to remote locations using homogeneous drones. In this study, charging stations are incorporated into the problem to increase the drones' fly ranges. They used a linear function dependent on the number of packages that are loaded on a drone for assessing the drone's energy consumption flying between two particular points. The problem was formulated as an MILP model. In all of the examples, they considered 1 depot, 1 recharging station, and 5 demand locations.

In contrast to the studies surveyed above where drones' energy consumption model is formulated as a linear function, Cheng et al. (2020) studied a DRP with a nonlinear convex function of the payload weight. The authors designed a B&C algorithm to solve the problem where the nonlinear constraints were relaxed first, and then logical cuts and sub-gradient cuts were added in the solution process to ensure that the drones' energy restrictions are satisfied. They introduced two new sets of benchmark instances with up to 50 customers and showed that in some instances, the linear approximation method could generate infeasible trips in the optimal solution, because of excess energy consumption. However, the energy spent during departure and landing phases of flight as well as the effect of other factors such as wind's speed and direction in the energy consumption rate were ignored in their study. Radzki et al. (2019) studied the drone energy consumption indirectly. They considered the influence of wind direction on the drones' energy consumption rate. They assumed that a drone's speed was constant during the flight, and the wind speed and direction remained unchanged not affecting on the drone's flight time. They also ignored the effect of a drone's payload weight on the energy consumption rate. The problem was implemented in a constraint programming environment and solved for examples with at most 9 customers.

As the accurate modeling of energy consumption makes the resulting model difficult to solve, the studies surveyed above tried to simplify the drone's energy consumption rate. However, with all the simplifications, the dimensions of the exactly solved DRP instances are still much smaller compared to similar VRP instances. Different energy models in DRPs lead to wide differences in a drone's fly range, and the routes obtained in the optimal solutions may be very different for the same instance. Furthermore, an energy consumption model may provide good approximations under specific conditions, but it generate infeasible routes due to excess energy consumption in reality, which consequently results in drone crashes and environment damages.

If one considers only the effect of payload weight on a drone's power consumption (and ignores the other factors and flight phases), the drone's energy consumption rate in the level flight (i.e., at a constant altitude and constant speed) can be formulated as a nonlinear convex function of the drone's total weight including the payload and battery weights (Stolaroff et al., 2018; Cheng et al., 2020; Kirschstein, 2020, Langelaan et al., 2017). In general, to better estimate the total energy consumption of a drone delivery mission, it is important to determine the energy consumed in all the phases of a delivery and return flights including takeoff, level flight, hovering, and landing considering other factors such as avionics, wind direction, and speed.

From a modeling perspective, the studies on DRPs only extended a known formulation of the Vehicle Routing problem (VRP) to fit special characteristics of drones. Although DRPs and VRPs have many common attributes, but drones' restrictions on flight time and payload capacity fundamentally change the solution space of the problem. The fly range of a drone is most dependent on its power consumption rate. Therefore, there is a need for a new formulation approach that fits better with drones' special characteristics; an approach that can incorporate more precise energy consumption models and can solve the problem efficiently.

## 2.2. Contributions of current paper

This paper contributes to the literature as follows:

(i) **General energy consumption model:** For the first time, this paper presents a new formulation approach for developing an MILP model for the DRP with a general energy consumption model after a smart preprocessing of the problem's input data. The new approach allows us to use the best available model, instead of assuming that flight range is a fixed value or using linear or convex approximations.

(ii) **Computational efficiency:** The MILP obtained by the new approach can be solved efficiently, while it is developed for a general energy function. It can solve the benchmark instances with linear and convex energy-consumption models on average 20 times faster than a recent B&C algorithm.

(iii) **Extendibility to challenging flight settings:** The MILP model can be simply extended to consider the energy consumed by the drones hovering at customer locations for a longer period when they arrive earlier before the opening of the time window. This is useful whenever the delivery company avoids early land due to safety issues or the customers do not allow drones to land out of the time window. Moreover, the linear formulation can be simply extended to cover a load-dependent flight policy in which the flight's controllable factors can depend on both path's characteristics and payload weight. In this case, both flight time and energy consumption amount for flying each arc may depend on the payload weight. It is the first time that such extensions are cast as linear formulations.

The rest of the paper shows how these contributions are achieved.

## 3. Problem statement and classical formulation

The drone routing problem (DRP) considered here is formally stated in Section 3.1. The notation needed to formulate the problem is given in Section 3.2. The classical formulation of the problem is presented in Section 3.3.

### 3.1. Problem statement

The DRP studied in this paper is a known problem studied by different authors with slight differences. In this problem, a set of homogeneous drones is considered. Each drone needs to start from a start depot and returns to an end depot in a limited time period while delivering parcels to customers on its route within known time windows. Each drone can travel for multiple trips where the total demand of the customers served in each trip is less than or equal to the drone's capacity. Drones are dispatched from the start depot with full energy batteries and new loads.

The drones have navigational abilities, and they can autonomously follow the commands and avoid collisions in the flight path. The estimated fly time between any two locations is a known given parameter. Without loss of generality, the time needed for swapping or recharging a drone's battery and for loading and unloading parcels on/off the drone is neglected as these times can be added simply to the flight time between any two points.

The objective of the problem is to minimize the total logistics cost that consists of the energy consumption cost and the transportation cost of drones. The problem determines the routes traveled by the drones and the time scheduling of each itinerary for serving each customer in a predetermined time window while satisfying drones' battery and payload limitations.

In this paper, the words *trip* and *route* are distinguished in a way that each trip starts from the start depot, continues visiting and serving customers and returns to the end depot; while the route of a drone consists of the sequence of the trips done by the drone. Moreover, the location of both start and end depots are considered the same; they are separated to simplify the modeling procedure.

The DRP described above is formulated as a mixed-integer nonlinear programming model in Section 3.3 (the old formulation) and mixed-integer linear programming (MILP) model in Section 4.2 (the new formulation).

### 3.2. Notation

To represent the delivery environment of the DRP, a directed graph $G_O = (V_O, A_O)$, called the *original* graph in the paper, where $V_O$ consists of three disjoint sets, i.e., $V = \{s\} \cup V_O^- \cup \{e\}$, is used. A vertex $i \in V_O^- = \{1,2,...,n\}$ signifies delivery location of (delivery) request $r_i$, and vertexes $s$ and $e$ are representative of the start depot and end depot locations, respectively. A fly time $t_a$ and a transportation cost $c_a$ is assigned to each arc. The rest of the notation is organized in the following subsections.

### 3.2.1. Sets and Parameters

The sets, indexes, and parameters used for the mathematical representation of the DRP are given below. Carefully note that indexes $v$ and $a$ are used for indexing vertexes and arcs in both original and *generated* graph (proposed in the next section), while the vertex and arc sets of both graphs differ.

| | |
|---|---|
| $n$ | The number of (delivery) requests |
| $s$ | The start depot |
| $e$ | A copy of $s$ that denotes the end depot, used to distinguish the start point and end point of each trip |
| $R$ | The set of requests, $R = \{r_1, r_2, \ldots, r_n\}$ |
| $V_O$ | The overall set of vertexes of the original graph, $V_O = V_O^- \cup \{s\} \cup \{e\}$ |
| $V_O^-$ | The overall set of vertexes of the original graph related to delivery locations of requests, i.e. $V_O^- = \{1, 2, \ldots, n\}$ |
| $A_O$ | The overall set of arcs of the original graph, i.e., $A_O \coloneqq \{(i,j): i \in \{s\}, j \in V_O^- \text{ or } i \in V_O^-, j \neq i \in V_O^- \cup \{e\}\}$ |
| $\delta^{in}(v)$ | The set of incoming arcs of vertex $v \in V_O$, i.e., $\delta^{in}(v) \coloneqq \{(w, u) \in A_O: u = v\}$ |
| $\delta^{out}(v)$ | The set of outgoing arcs of vertex $v \in V_O$, i.e., $\delta^{out}(v) \coloneqq \{(w, u) \in A_O: w = v\}$ |
| $N$ | The number of available drones |
| $q_r$ | The demand weight of request $r \in R$ |
| $r(v)$ | The request associated with vertex $v \in V_O$ |
| $Q$ | The capacity of each drone in kg |
| $M$ | The battery capacity of each drone in units |
| $t_a$ | The total fly time of a drone on passing arc $a \in A_O$ in minutes |
| $c_a$ | The logistic cost of passing arc $a \in A_O$ in unit |
| $\delta$ | The energy cost in unit |
| $a_r$ | The earliest time that request $r \in R$ can be delivered |
| $b_r$ | The latest time that request $r \in R$ can be delivered |
| $a_d$ | The earliest permitted departure time of drones from the starting depot |
| $b_d$ | The latest allowed arrival time of drones to the end depot. |

### 3.2.2. Decision and auxiliary variables

To formulate the DRP, six different sets of variables are defined as follows:

$x_a$      A binary decision variable that is equal to 1, if a drone flies arc $a \in A_O$; and 0 otherwise

$w_a$      A nonnegative variable that specifies the payload of the drone in flying arc $a \in A_O$

$ed_a$      Energy depletion of a drone's battery after flying arc $a \in A_O$

$f_v$      A nonnegative variable representing the energy consumed from the battery of a drone at the start of serving vertex $v \in V_O$

$y_v$      A nonnegative variable indicating the time at which a drone starts service at vertex $v \in V_O$

$z_{vw}$      A binary variable that is equal to 1 if the same drone performs the trip whose first request is $r(w)$ immediately after the trip whose last request is $r(v)$, and 0 otherwise; $v, w \in V_O^-$.

### 3.3. Old formulation

Using the notation introduced in Section 3.2, the classical formulation of the DRP is as follows:

$$\min \sum_{a \in A_O} c_a x_a + \delta \sum_{a \in A_O} ed_a x_a \tag{1}$$

subject to:

$$\sum_{a \in \delta^{in}(v)} x_a = 1 \qquad \forall v \in V_O^- \tag{2}$$

$$\sum_{a \in \delta^{out}(v)} x_a = 1 \qquad v \in V_O^- \tag{3}$$

$$\sum_{a \in \delta^{out}(s)} x_a = \sum_{a \in \delta^{in}(e)} x_a \tag{4}$$

$$\sum_{a \in \delta^{out}(v)} w_a = \sum_{b \in \delta^{in}(v)} w_b - q_{r(v)} \qquad v \in V_O^- \tag{5}$$

$$w_a \leq Q x_a \qquad a \in A_O \tag{6}$$

$$\sum_{a \in \delta^{in}(e)} w_a = 0 \tag{7}$$

$$f_w - ed_a \geq f_v - M_a^1(1 - x_a) \qquad a = (v, w) \in A_O \tag{8}$$

$$f_e \leq M \tag{9}$$

$$f_s = 0 \tag{10}$$

$$y_v + t_a \leq y_w + M_a^2(1 - x_a) \qquad a = (v, w) \in A_O \tag{11}$$

$$y_v + t_{a=(v,e)} + t_{a=(s,w)} \leq y_w + M_{vw}^3(1 - z_{vw}) \qquad v \neq w \in V_O^- \tag{12}$$

$$a_{r(v)} \leq y_v \leq b_{r(v)} \qquad v \in V_O^- \tag{13}$$

$$a_d \leq y_v \leq b_d \qquad v \in \{s, e\} \tag{14}$$

$$\sum_{v \neq w \in V_O^-} z_{vw} \leq x_{a=(s,w)} \qquad w \in V_O^- \tag{15}$$

$$\sum_{w \neq v \in V_O^-} z_{vw} \leq x_{a=(v,e)} \qquad v \in V_O^- \tag{16}$$

$$\sum_{a \in \delta^{out}(s)} x_a - \sum_{v \in V_O^-} \sum_{w \neq v \in V_O^-} z_{vw} \leq N \tag{17}$$

$$\sum_{a \in \delta^{out}(s)} x_a \geq \frac{\sum_{i \in R} q_i}{Q} \tag{18}$$

$$x_a \in \{0,1\} \qquad a \in A_O \tag{19}$$

$$z_{vw} \in \{0,1\} \qquad v \neq w \in V_O^- \tag{20}$$

$$ed_a = F_a(w_a, t_a) \qquad a \in A_O \tag{21}$$

$$f_v \geq 0 \qquad v \in V_O \tag{22}$$

$$y_v \geq 0 \qquad v \in V_O \tag{23}$$

$$w_a \geq 0 \qquad a \in A_O \tag{24}$$

where $M_a^1$, $M_a^2$, and $M_{vw}^3$ are sufficiently large positive numbers.

The first and second terms in the objective function (1) aim to minimize the transportation and energy depletion cost of passing arcs for delivering requests. Constraints (2) guarantee that each request is delivered exactly once. Constraints (3) ensure the flow conservation at the vertexes and make sure that the out-degree of a vertex is equal to its in-degree. Constraints (4) indicate that the number of trips leaving the start depot is equal to the number of arriving drones at the end depot. Constraints (5) track the payload of the drones on their fly paths, and constraints (6) ensure that the payload capacities of drones are not exceeded. Constraints (7) make sure that the drones return to the end depot with empty load. Constraints (8) enable the model to track the drones' consumed energy during the delivery operations. Constraints (9) limit the energy consumed to the battery capacity of the drones. Constraints (10) set the consumed energy at the start of each trip to zero. Constraints (11) enable the model to track the time during the delivery operations. Constraints (12) set the service time relationship between two successive trips performed by the same drone. Constraints (13) and (14) describe the time window restrictions. Constraints (15) and (16) define the relation between variables $x$ and $z$. Constraints (17) limit the number of drones that are available. Constraint (18) is a valid inequality specifying the minimum number of drones needed for delivering the parcels. Finally, constraints (22)–(24) are domain constraints.

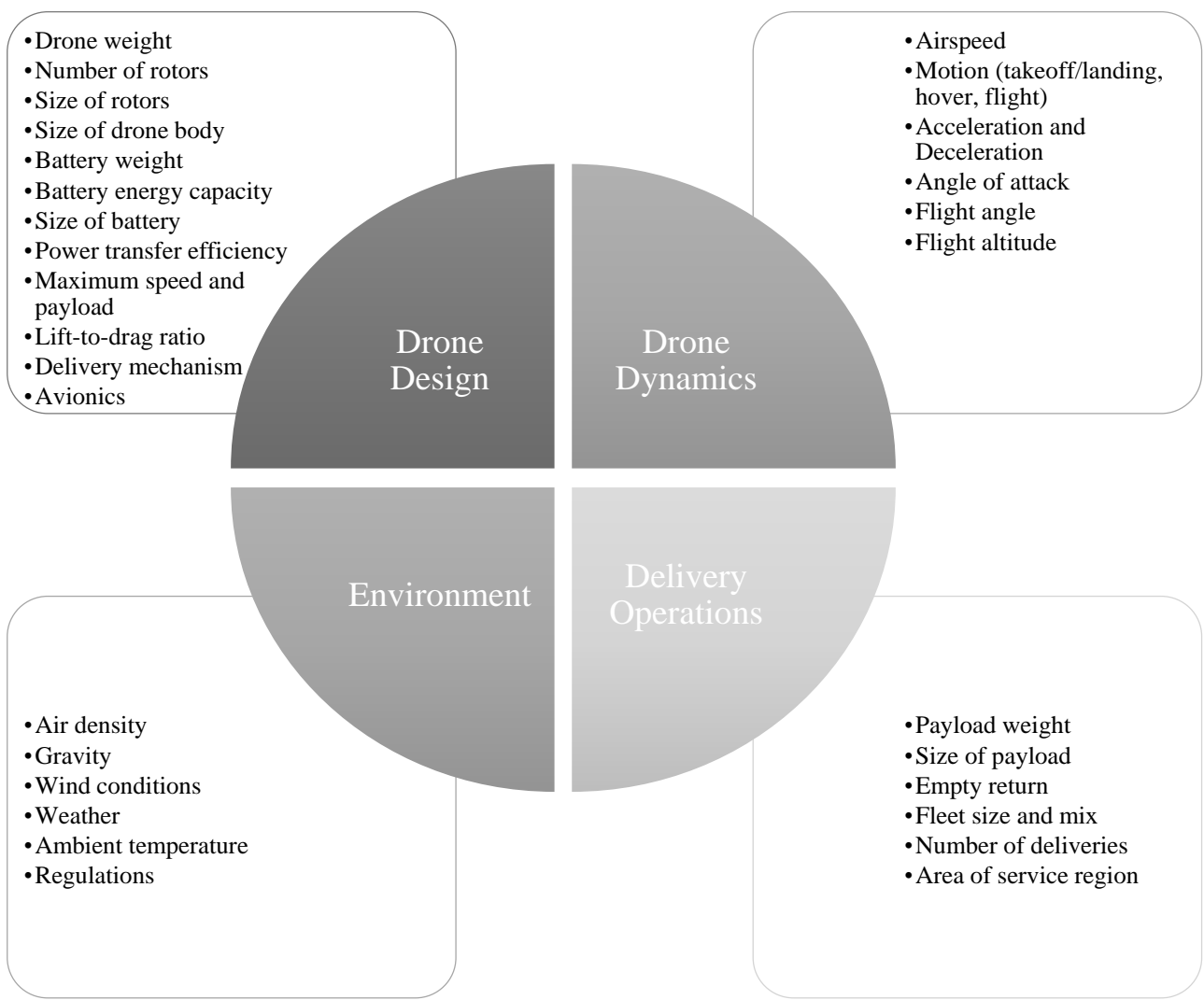

**Figure 1.** Factors effecting a drone's energy consumption

Constraints (21) consider the interaction of a drone's variable weight and its energy consumption. One can see this set of constraints makes the formulation nonlinear.

The function $F_a$ in (21) calculates a drone's consumed energy for flying arc $a \in A$, which is usually determined based on different simplifying assumptions. Common convex and linear functions that approximate energy consumption are presented in equations (25) and (26), respectively,

$$F_a(w_a, t_a) = t_a \times \alpha(\beta + w_a)^{3/2} \tag{25}$$

$$F_a(w_a, t_a) = t_a \times (\partial w_a + \theta) \tag{26}$$

where $\alpha, \beta, \partial$, and $\theta$ are given parameters.

The above formulation has been used by different papers with small differences. For example, Dorling et al. (2017) used a similar model whose objective is minimizing the makespan ($\min y_e$) where the energy consumption function is an affine function of $w_a$, given in (26). Cheng et al. (2020) considered a similar formulation and proposed a B&C method for solving the above model when $F_a$ is a convex function of $w_a$, which can be any of (25) and (26).

Zhang et al. (2020) classified all the influencing factors of the energy consumed during a drone's flight into four groups represented (see Figure 1). Actually, the energy consumption models given in (25) and (26) are very naive and ignore many factors compared to the best available choices, which have been continually improved over the years. For example, Kirschstein (2020) studied the following model (initially proposed in Langelaan et al. (2017))

$$F_a(w_a, t_a) = \sum_{i=1}^{4} t_i G_a(w_a, i) \tag{27}$$

with $t_a = \sum_{i=1}^{4} t_i$. The index $i = 1, \ldots, 4$ indicates each phase (i.e., takeoff, level, hovering, and landing, respectively) of the idealized delivery flight of a drone flying arc $a$ (note that for a return flight to end depot the same formula can be used by setting $w_a = 0$). The time required for completing phase $i$ of the flight is denoted by $t_i$, which is determined based on the flight policy (flight angle, flight velocity, and distance covered in phase $i$). For each phase $i$ of the flight, the function $G_a(w_a, i)$ is defined as follows (here an error in the original formula has been fixed after contacting the authors of Langelaan et al., (2017)):

$$G_a(w_a, i) = \alpha_1^i + \alpha_2^i \phi^i[w_a]\sqrt{\alpha_3^i + \alpha_4^i(\beta + w_a)^2} + \alpha_5^i(\beta + w_a)^{3/2} + \alpha_6^i(\beta + w_a)^{1/2} + \alpha_7^i(\beta + w_a)$$

where $\phi^i[w_a]$ is obtained by solving the following $w_a$-dependent equation:

$$\alpha_8^i \sqrt{\alpha_3^i + \alpha_4^i(\beta + w_a)^2} = X \sqrt{\alpha_9^i + X^2 + X \frac{(\alpha_{10}^i + \alpha_{11}^i(\beta + w_a))}{\sqrt{\alpha_{12}^i(\beta + w_a)^2 + \alpha_{13}^i + \alpha_{14}^i(\beta + w_a)}}}$$

with respect to $X$. The parameters $\alpha_j^i$, $j = 1, ...,14$, are independent of payload weight $w_a$ and determined based on the *environmental factors* (gravity, wind speed and direction, air density, weather conditions), drone *design specifications* (drone weight, number of rotors, size of rotors, rotor mean chord, number of blades, blade lift, size of drone body, battery weight, size of battery, rotor solidity ratio, power transfer efficiency), and other pre-set *controllable factors* (the flight angle, altitude, and speed at each phase $i$). In Appendix A, the parameters $\alpha_j^i$, $j = 1, ...,14$, are rewritten in terms of the parameters originally used by Kirschstein (2020); here they are used instead of the original ones to simplify the presentation.

One can see that the energy-consumption model in (27) is a very complex non-convex function of payload weight $w_a$, which cannot be handled by the methods developed for (25) or (26); in fact, (25) becomes functionally equivalent to (27) if one ignores all the other terms, i.e., by setting $\alpha_1^i = \alpha_2^i = \alpha_6^i = \alpha_7^i = 0$, $i = 1, ...,4$. Fortunately, our solution method alleviates the way of using such an energy consumption model.

## 4. New formulation of drone routing problem

This section provides the details of our new formulation. Section 4.1 first presents a procedure to construct a new graph that is used as an input for the new formulation, and then Section 4.2 proposes the MILP formulation.

### 4.1 Graph generation and preprocessing phase

The new formulation is defined on a directed incomplete graph $G = (V, A)$, called the *generated* graph hereafter, where $V$ is the vertex set and $A$ is the arc set. The vertexes of the graph $G_O = (V_O, A_O)$ used in the old formulation of the DRP, represent only the geographic information of the locations where a drone can land or dispatch to deliver the requests. However, besides the geographic positions, the vertexes of the graph $G = (V, A)$ include data about the drones' payloads when arriving at these locations. To construct the new graph $G = (V, A)$, many copies of each delivery vertex in the original graph $G_O = (V_O, A_O)$ are created, where in each copy a possible set of requests loaded on a drone **reaching** the delivery location is considered. Hence, a delivery vertex $v \in V^-$ of the graph $G$ can be represented by $v = (r, L)$ with $r \in R$ where $L \subseteq R$ is any subset of the requests that includes request $r$ ($L$ must contain $r$ because of the problem assumption that allows a drone to visit any customer location only once). Note that if a drone reaches vertex $v = (r, L) \in V^-$, the drone is located in the geographic location of request $r$ while its payload consists of the requests included in the set $L$.

To simplify the presentation, the set of vertexes that are related to the request $r \in R$ is represented as *cluster $V_r$* defined below:

$$V_r := \{(r, L) : L \subseteq R \text{ with } L \ni r \text{ and } feas_{r,L} = 1\} \tag{28}$$

where $feas_{r,L}$ is a 0-1 parameter that equals 0 if the vertex $(r, L)$ is recognized as infeasible. For each of the start and end depots, the corresponding clusters are defined as $V_s = \{s\}$ and $V_e = \{e\}$, respectively. Therefore, the whole set $V$ of vertexes is given by $V = \cup_{r \in R} V_r \cup \{s\} \cup \{e\}$, and $V^- = \cup_{r \in R} V_r$.

In (28), $feas_{r,L}$ must be initially set to 1 for any $L \ni r$, but by applying appropriate rules one can reduce the size of the generated graph by setting $feas_{r,L} = 0$ for some sets $L \ni r$. Three rules used in this paper are listed below (see Appendix B for additional technical details):

R1. Full checking of capacity constraint: the drone's capacity limitation <u>is satisfied</u> for $L$, that is, the sum of the weights of requests in the set $L$ must satisfy the drone's capacity constraint (here we can also check other constraints on the number of carried requests, shape, etc.).
$feas_{r,L} = 0$ if the requests in set $L$ does not satisfy the physical capacity of the drone.

R2. Partial checking of battery and time-window constraints: <u>for at least one permutation</u> of the requests contained in the set $L$, serving the requests in the order of the permutation is feasible in a single trip <u>assuming</u> that the drone carrying all the requests in $L$ flies from the depot directly to the current vertex, that is, the drone battery constraint and customer time-window constraints are all satisfied in the rest of the flight for the permutation considering direct fly from the depot to the current position (which is necessary condition for that the battery and time-window constraints are satisfied).
$feas_{r,L} = 0$ if no permutation satisfying the above condition can be found.

R3. If $feas_{r,L}$ becomes 0 for some $L$ (due to rule R1 or R2), then $feas_{r,L'} = 0$ for all the <u>super</u>sets $L' \supseteq L$.

Carefully note that one can define all possible vertices for any set $L \ni r$ that only satisfies the drones' capacity limitations (R1 and R3) <u>without</u> partially checking the drones' battery and customers time-window constraints based on the second rule R2; however, in this case, the generated graph becomes very large. Also, remark that the <u>full checking</u> of the battery and time-window constraints is done in the main formulation given in Section 4.2.

The arc set $A$ of the generated $G$ is defined as the set of any possible moves between pairs of vertexes in $V$. Actually, a directed arc $a \in A$ is only generated between two vertexes, if a drone moves from one vertex to another one, the change in the drone's payload makes sense. Therefore, the arc set A is defined as follows:

I. Arcs $(a,b) \in V \times V$ outgoing from any vertex $a = (r,L) \in V^-$ with $|L| > 1$ can enter to a vertex $b = (p,K) \in V^-$ with $K = L \setminus \{r\}$.

II. Arcs $(a,b) \in V \times V_e$ outgoing from any vertex $a = (r,L) \in V^-$ with $|L| = 1$ can enter to end depot $b \in \{e\}$.

III. Arcs $(a,b) \in V_s \times V$ outgoing from the start depot $a = s$ can enter to any vertex $b = (r,L) \in V^-$.

As mentioned earlier, drones have severe constraints in payload capacity and maximum flight duration due to battery depletion in practice, which causes drones to be capable of loading a small number of packages in each single trip. Therefore, for many subsets $L$ of the set $R$, $feas_{r,L}$ can be set to 0. Hence, in practice, the number of the vertexes does not grow exponentially, and thus such a graph construction makes sense when using drones for delivering parcels.

The preprocessing algorithm for constructing the generated graph and its main sub-algorithms are presented in the next three subsections.

### 4.1.1. The preprocessing algorithm for generating the new graph

The pseudocode of the preprocessing algorithm used to generate vertexes and arcs of the graph used in the new formulation based on the guidelines discussed earlier is given below.

**The preprocessing algorithm for constructing the generated graph**

1) **Input** $R, s, e, Q, M, [a_d, b_d], [a_r, b_r], q_r$ for $r \in R$;
2) Set $FeasibleLoads = \emptyset$;
3) Set $Load = \emptyset$;
4) Set $choices\_left = |R|$;
5) Set $index = 0$;
6) $FeasibleLoad$=FindSubset ($R, hoices\_left, index, Load, FeasibleLoad, s, e, Q, [a_d, b_d], [a_r, b_r], q_r$ for $r \in R$ );
7) **For** $r \in R$, set $V_r = \emptyset$;
8) **For** each set $L \in FeasibleLoad$
9)     For each request $r \in L$
10)        Generate a vertex $v = (r,L)$;
11)        $V_r = V_r \cup \{v\}$;
12)    **End For**
13) **End For**
14) Set $V^- = \cup_{r \in R} V_r$ and $V = \cup_{r \in R} V_r \cup \{s\} \cup \{e\}$;
15) $A = \emptyset$;
16) **For** each $v = (r,L) \in V^-$
17)    Generate arc $a = (s,v)$;
18)    Set $A = A \cup \{a\}$;
19)    **If** $|L| == 1$;
20)        Generate arc $a = (v,e)$ ;
21)        Set $A = A \cup \{a\}$;
22)    **End If**
23)    **For** each $w = (p,I) \in V^-$
24)        **If** $I = L \setminus \{r\}$
25)            Generate arc $a = (v,w)$;
26)            Set $A = A \cup \{a\}$;
27)        **End If**
28)    **End For**

29) **End For**
30) **Return** $G = (V, A)$

The first steps of the algorithm initiate sets, parameters, and variables. In Step 6 using the FindSubset algorithm (outlined in Section 4.1.2), all the sets $L \subseteq R$ with $feas_{r,L} = 1$ are generated and saved in $FeasibleLoad$. Next, in Steps 8–13, for each set $L \in FeasibleLoad$ and for each request $r \in L$, a vertex is created and included in the cluster $V_r$. After the vertex set $V$ of the graph is generated, in Steps 16–18, the arcs of the new graph are generated that start from the vertex related to the start depot to any vertex in the graph except the vertex associated with the end depot. Steps 19–22 define arcs starting from any vertex $v \in V$ with $|L| = 1$ to the vertex associated with the end depot. Eventually, in Steps 23–28, an arc is created from any vertex $v = (r, L) \in V^-$ to vertex $(p, I) \in V^-$ if and only if $I = L \setminus \{r\}$. The algorithm ends by returning the *generated* graph $G = (V, A)$.

### 4.1.2. Algorithm for fining any set of requests that can be served in a single trip

The FindSubset algorithm used in the preprocessing algorithm is a recursive function that uses its own previous term to generate any feasible subset $L$ of the set $R$ that can be serviced in a single trip. It uses the function IsLoadPossibble (outlined in Section 4.1.3) to determine whether serving customers contained in a given set is possible in a single trip or not. The algorithm works considering that if delivering a set of requests $L \subseteq R$ is not possible in a single trip, then delivering any super set of $L$ is also impossible (rule R3). The pseudocode of the FindSubset algorithm is represented in the following.

**The FindSubset algorithm**

1) **Input** $R, choices\_left, index, Load, FeasibleLoad, s, e, Q, M, [a_d, b_d], [a_r, b_r], q_r$ for $r \in R$ ;
2) **If** IsLoadPossibble( $Load, s, e, Q, M, [a_d, b_d]$ and $[a_r, b_r]$ and $q_r$ for $r \in Load$)
3)     **If** ($|Load| > 0$)
4)         Set $FeasibleLoad = Load \cup FeasibleLoad$;
5)     **End If**
6)     **If** ($choices\_left > 0$)
7)         **For** $i = index + 1$ to $i = |R|$
8)             $Load = Load \cup \{r_i\}$;
9)             Set $Load2 = Load$;
10)            FindSubset ($R, choices\_left - 1, i + 1, Load, FeasibleLoad, s, e, Q, M, [a_d, b_d], [a_r, b_r], q_r$ for $r \in R$);
11)            Set $Load = Load2$;
12)         **End For**
13)     **End If**
14) **End If**
15) **Return** $FeasibleLoad$

### 4.1.3. Algorithm for determining whether a given set of requests can be served in a single trip

The algorithm IsLoadPossibble checks the sum of the weights of the requests included in the set $Load$ and returns false if the weight exceeds drones' maximum capacity. Next, it tries to find a permutation of the requests contained in the set $Load$ where delivering requests in the sequence of the permutation *partially* satisfies drones' limit in consumed energy and customers' time windows (rule R2). The algorithm ends by returning true if a feasible permutation is found, and false otherwise. In this algorithm, ConsumedEnergy is a black-box model that returns the energy needed for delivering requests in the predetermined order; it can work based on an explicit formulation (e.g., (25)), a computational procedure (e.g., (27)), or a simulation model. The pseudocode of the IsLoadPossibble algorithm for determining whether a given set of requests can be served in a single trip or not is given below.

---

**The IsLoadPossibble algorithm**

1) **Input** $Load, s, e, Q, M, [a_r, b_r], q_r$ for $r \in Load$;
2) Set Temp = 0;
3) **For** each $r \in Load$
4)     $Temp = Temp + q_r$;
5) **End For**
6) **If** $Temp > Q$
7)     **Return False**;
8) **End If**
9) Set Found = False;
10) **For** each permutation $per$ of the set $Load$ do:
11)     **If** ConsumedEnergy($per, s, e,$ and $q_r$ for $r \in Load2) \leq M$
12)         **If** serving each requests $r \in Load$ is possible in in $[a_r, b_r]$ considering the sequence listed in $per$
13)             **Set** Found=True;
14)             **Return** True;
15)         **End If**
16)     **End If**
17) **End For**
18) **If** Found == False
19) **Return** False;
20) **End If**

---

### 4.1.4. Pre-analysis of the generated graph's size

In the original graph $G_O = (V_O, A_O)$, for each request one vertex exists; however, in the *generated* graph $G = (V, A)$, for each location several vertexes are created. Certainly, this enlarges the new graph, and consequently the number of variables and constraints of the new linear formulation (proposed in Section 4.2) become greater compared to the old nonlinear formulation (Section 3.3).

In the following, some reasonable bounds on the number of vertexs and arcs of the generated graph are derived in (29) and (30), respectively, in a pre-analysis stage (that is, before constructing the generated graph). This enables us to compare them with the actual ones obtained after finishing the preprocessing phase (post-analysis stage), which reveals why the new approach is efficient in practice. In fact, this shows the perfect performance of applying the rules R1, r2, and R3 in reducing the size of the generated graph.

Let $UPMNR$ denote an easy-to compute upper bound on the maximum number of requests that each drone can serve in a single trip (denoted by $MNR$). Computing $MNR$ for a given instance is not easy, hence $UPMNR$ is used in our pre-analysis (see Appendix C for more details). Then, an upper bound on the number of vertexes of the generated graph can be calculated as follows:

$$|V| \leq 2 + |R| \sum_{i=0}^{UPMNR-1} \binom{|R|-1}{i}. \tag{29}$$

It is clear that $G$ is an *incomplete directed* graph where the number of arcs out-going from any vertex $v = (r, L) \in V^-$ with $|L| > 1$ is at most $|L| - 1$. Therefore, one can consider the following upper bound on the cardinality of the arc set of the generated graph $G = (V, A)$:

$$|A| \leq |R| \sum_{i=0}^{UPMNR-1} \binom{|R|-1}{i} + |R| \left( \sum_{i=0}^{UPMNR-1} \binom{|R|-1}{i} \times i \right) + |R|. \tag{30}$$

According to the above equations, the number of vertexes and the number of arcs are bounded by $O(|R|^{UPMNR})$ and $O(UPMNR \times |R|^{UPMNR})$, respectively, where $UPMNR$ can be $|R|$ in the worst case. These bounds are meaningfully greater than the number of the vertexes and the number of arcs of the original graph used in the classical model, which are bounded by $O(|R|)$ and $O(|R|^2)$, respectively. However, it is not a problematic issue in practice as our extensive numerical study shows this in Section 5.1. Indeed, when using drones for parcel delivery, usually for many subsets $L$ of $R$, $feas_{r,L}$ becomes 0, and consequently the actual dimensions of the generated graph are much smaller than the bounds obtained in the worst case.

In Section 5.1, the size of the generated graph for the benchmark instances is evaluated and compared with the upper bounds obtained using equations (29) and (30), which show only small percentages of these bounds are needed in practice. Moreover, the times required for preprocessing phase is fairly small compared to the solver times (see Figure 4 for a summary).

### 4.1.5. A simple example of the generated graph

In the generated graph $G = (V, A)$, the vertex set $V$ represents the geographical locations and the drones' payloads when they visit these locations. A feasible solution to the problem consists of at most $N$ drone paths where each path consists of one or more trips that start from the start depot and ends at the end depot. Figure 2 depicts an example of the generated graph when $R = \{1, 2, 3, 4\}$, $q_1 = q_2 = q_3 = 0.3$, $q_4 = 0.2$, $Q = 0.5$, and $[a_i, b_i] = [0, +\infty)$ for $i \in \{1,2,3,4\}$. In this figure, each vertex is represented by a square and different clusters are specified using dashed circles. As one can see, a drone can move from the start depot to any vertex, and only one arc is defined from each cluster to the end depot. Vertexes with the same color are those that the drone reach them with the same drone payload.

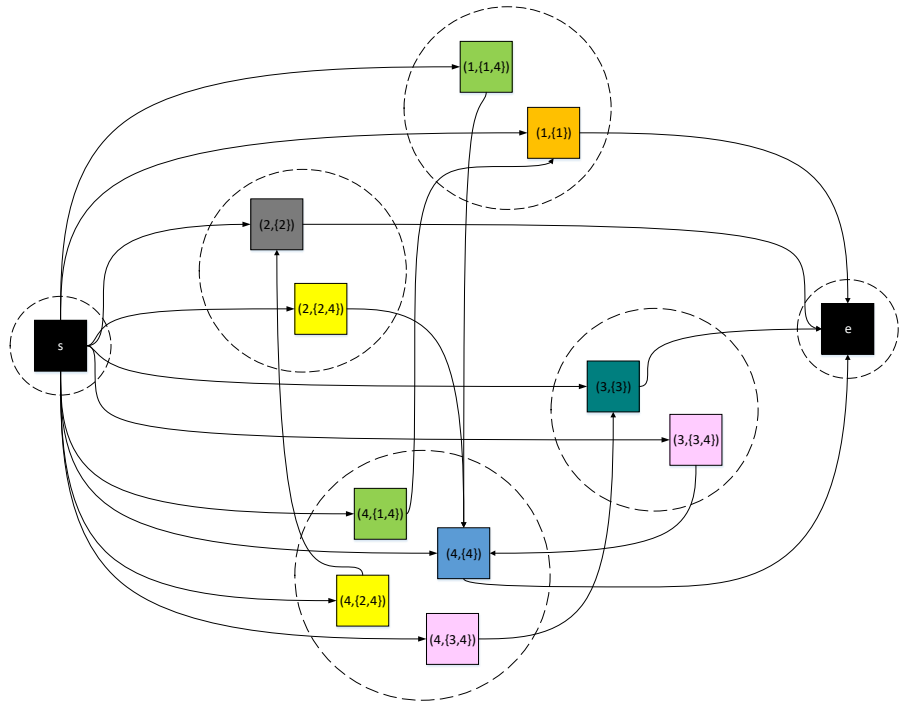

**Figure 2.** The generated graph for an instance with $|R| = 4$, $q_1 = q_2 = q_3 = 0.3$, $q_4 = 0.2$, and $Q = 0.5$.

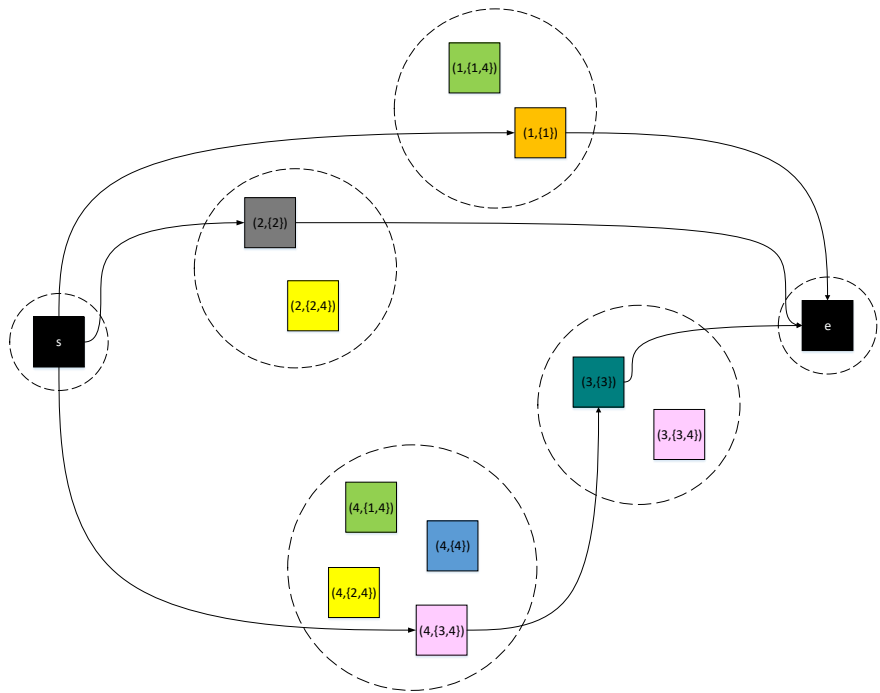

**Figure 3.** A feasible solution for the generated graph given in Figure 2.

Figure 3 displays a selected set of arcs that form a feasible solution for the generated graph given in Figure 2. For each feasible solution, one arc enters to each cluster and outgoes from it.

### 4.1.6. Practical advantages of using generated graph

The main challenge of using drones for delivery is how to ensure on time and failure-free package delivery, which largely depends on accurately calculating the drones' energy consumption. Therefore, the designed graph has the benefit that perusing the energy of a drone moving in the graph is more precise and much easier as well. For example, one can evaluate the energies needed in different drone's motions such as takeoff, landing, hover, and horizontal flight separately, and then consider the sum as the energy required for passing the arc (see (27) for example). Moreover, the effect of the weight and size of the payloads and batteries, different weather conditions and wind directions (whenever their changes are not significant during the delivery), as well as empty returns can be handled in the generated graph. In fact, using the generated graph, the mathematical model becomes independent of the model used for calculating the energy consumed when moving from any location to another one in the optimization phase (though it indirectly depends on the energy consumption model via the calculated parameters). Moreover, recall that the payload capacity constraints, which must be explicitly formulated in the classical model, are implicitly managed in the generated graph.

### 4.2 New formulation

This subsection formulates the DRP stated in Section 3 as an MILP model using the generated graph described in Section 4.1. The problem statement remains unchanged *excepting* that the energy-consumption computing procedure can be chosen arbitrarily.

Some new notation and modifications are required to present the new formulation. Carefully **note** that here, with a little abuse of notation, indexes $v \in V$ and $a \in A$ are similarly used to indicate a vertex or an edge in the generated graph $G = (V, A)$, which should not be confused with $v \in V_O$ and $a \in A_O$ used in Section 3.3.

As the generated graph is used in the new model, the variables $ed_a$ appeared in the old formulation (representing the energy depletion of a drone's battery after flying arc $a \in A$) are replaced by the constants because the payload of the drone flying each arc is known in the generated graph at each vertex; in other words, each $ed_a$ is a parameter in the new formulation, and the nonlinearity caused by $ed_a$ through (21) is removed.

Let $cl(v)$ represent the cluster number associated with vertex $v \in V$, defined as follows

$$cl(v) = \begin{cases} i & v \in V_{r_i} \\ 0 & v = s \\ n+1 & v = e \end{cases}.$$

The variables $w_a$ and $ed_a$ are not needed in the new formulation. The following changes are also required:

$x_a$    A binary variable that equals 1 if a drone flies arc $a \in A$, and 0 otherwise

$f_i$     The accumulated energy consumption of a drone upon arrival at **cluster** $i$; $i \in \{0, \ldots, n+1\}$

$y_i$     The time at which a drone starts providing service at **cluster** $i$ in minutes; $i \in \{0, \ldots, n+1\}$

$z_{ij}$     A binary variable that is equal to 1 if the same drone performs the trip whose first request is in cluster $j$ immediately after the trip whose last request is in cluster $i$, and 0 otherwise, $i \neq j \in \{1, 2, \ldots, n\}$.

Note that the variables $z_{ij}$ implicitly determine the route (which may be a sequence of multiple trips) that must be done by the same drone.

Using the notation introduced above, the DRP can be formulated as follows:

$$\min \sum_{a \in A} c_a \times x_a + \delta \sum_{a \in A} ed_a x_a \tag{31}$$

subject to:

$$\sum_{a \in \delta^{in}(v)} \sum_{v \in V_r} x_a = 1 \quad r \in R \tag{32}$$

$$\sum_{a \in \delta^{in}(v)} x_a = \sum_{a \in \delta^{out}(v)} x_a \quad v \in V^- \tag{33}$$

$$\sum_{a \in \delta^{out}(s)} x_a = \sum_{a \in \delta^{in}(e)} x_a \tag{34}$$

$$f_{cl(w)} - ed_a \geq f_{cl(v)} - M_a^4(1 - x_a) \quad a = (v, w) \in A \tag{35}$$

$$f_{n+1} \leq M \tag{36}$$

$$f_0 = 0 \tag{37}$$

$$y_i + t_{ij} \leq y_j + M_{ij}^5 \left(1 - \sum_{a=(v,w) \in A: cl(v)=i, cl(w)=j} x_a\right) \quad i \neq j \in \{0, 1, \ldots, n\} \tag{38}$$

$$y_i + t_{i,n+1} + t_{j,0} \leq y_j + M_{ij}^6(1 - z_{ij}) \quad i \neq j \in \{1, 2, \ldots, n\} \tag{39}$$

$$a_i \leq y_i \leq b_i \quad i \in \{1, \ldots, n\} \tag{40}$$

$$a_d \leq y_i \leq b_d \quad i \in \{0, n+1\} \tag{41}$$

$$\sum_{i=1:i\neq j}^{n} z_{ij} \leq \sum_{a=(v,w):\, v=s, w\in V_j} x_a \qquad j \in \{1,2,\ldots,n\} \tag{42}$$

$$\sum_{j=1:j\neq i}^{n} z_{ij} \leq \sum_{a=(v,w):\, v\in V_i, w=e} x_a \qquad i \in \{1,2,\ldots,n\} \tag{43}$$

$$\sum_{a\in\delta^{out}(s)} x_a - \sum_{i=1}^{n} \sum_{j=1:j\neq i}^{n} z_{ij} \leq N \tag{44}$$

$$\sum_{a=(i,j):\, i\in\{s\}, j\in V^-} x_a \geq \frac{\sum_{i\in R} q_i}{Q} \tag{45}$$

$$x_a \in \{0,1\} \qquad a \in A \tag{46}$$

$$z_{ij} \in \{0,1\} \qquad i \neq j \in \{1,2,\ldots,n\} \tag{47}$$

$$y_i \geq 0 \qquad i \in \{0,\ldots,n+1\} \tag{48}$$

$$f_i \geq 0 \qquad i \in \{0,\ldots,n+1\} \tag{49}$$

where $t_{ij}$ represents the flight time of a drone passing from the location of cluster $i$ to the location of cluster $j$ (recall that the vertexes in a cluster in the generated graph have the same location).

The way that the vertexes of the generated graph are defined, enables us to formulate the problem such that flying each arc not only determines the path of the drone, but also its payload when flying each arc on its path.

Objective function (31) minimizes the total transportation and energy depletion cost of passing the arcs by the drones for delivering the requests. Note that values of $ed_a$ are determined in the preprocessing stage. Considering the definition of clusters, constraints (32) guarantee that each customer is visited exactly once by a drone with an admissible payload that satisfies drone's maximum capacity limit. Constraints (33) ensure the flow conservation at the vertexes and make sure that the out-degree of any vertex is equal to its in-degree. These constraints also enable the model to track the payload of drones on their fly path. Constraints (34) indicate that the number of trips leaving the starting depot is equal to the number of trips arriving at the end depot. Constraints (35) enable the model compute the energy consumed by each drone during delivery operations. In these constraints, unlike the same constraints used in the old formulation, $ed_a$ is a parameter, and thus the nonlinearity that appears in (21) is removed in the new model. The descriptions of the next constraints are the same as the ones provided in Section 3.3, considering how the generated graph is created (Section 4.1). Recall that payload capacity constraints and constraints (21) (which explicitly use the energy consumption model) are not needed here as they are implicitly considered in the generated graph. Moreover, payload\arc-dependent battery constraints are replaced with only arc-dependent ones (35).

The big constants used in the above formulation are set as follows:

$$M_a^4 = M + ed_a$$

$$M_{ij}^5 = \max\{0, b_{r_i} - a_{r_j} + t_{ij}\}$$

$$M_{ij}^6 = \max\{0, b_{r_i} - a_{r_j} + t_{i,n+1} + t_{0,j}\}.$$

Note that with the new graph, constraint set (32) can be replaced by the following alternative:

$$\sum_{a=(v,w)\in A: v=s, w=(p,L): r\subseteq L} x_a = 1 \qquad r \in R. \tag{50}$$

Let us provide more explanations on how the new model consider *payload-capacity*, *battery-capacity*, and *time-window* constraints. The vertices in the generated graph keep the additional information on the set of requests loaded on the drone when arriving at a customer location. This information helps us to find out whether payload-capacity constraints are satisfied or not. In the generated graph, for calculating energy consumed flying an arc (not a trip), one just needs to know the payload weight of the drone flying the arc, which is known at the start vertex of the arc (knowing the order of the visited customers reaching to the location of the start vertex is not required). After calculating the energy required for flying each arc, in the new formulation, the total energy of a trip is controlled by linear constraints (35)–(37). In fact, the energy consumed flying each arc is calculated in the preprocessing phase whereas the energy consumed in each trip is determined in the new formulation, which depends on the order of the customers visited in the trip.

In order to fully check battery-capacity and time-window constraints, the order of customers to be visited in a trip has to be known. The above model control time-window constraints in (39)–(41). In the preprocessing phase, rule R2 is used to *partially* check battery-capacity and time-window constraints to decrease the size of the generated graph.

One can see that the new formulation always has a linear structure regardless of which type of energy consumption model is used. This property allows us to use powerful mixed-integer programming solvers for solving the MILP model efficiently considering the fact that the size of the generated graph does not exponentially grow compared to the original graph in practice. Section 5 assesses the computational efficiency of the new formulation.

## 4.3 Extension for considering extra energy consumed in early arrivals

To derive the old or the new formulations, it is assumed that the energy consumed while drones are on the ground is negligible, and that drones *are allowed* by both owner and customers to land and wait at customer sites whenever they arrive sooner. The formulation can be *simply* extended to include the energy consumed while drones wait on the ground (i.e., drones are alert even when landed on ground and ready for their next mission, and they use advanced technologies to detect dangerous situations such as stealing); it is easy because the energy used by a drone is independent of its payload's weight in this case.

Furthermore, the new formulation can be extended to another MILP in order to cover the more challenging case where the drones have to hover at customer locations if they arrive sooner due to safety issues or impossibility of landing out of the time windows; it is difficult because the energy used by a drone for a longer hovering period depends on its payload weight (note that extending the old formulation for this case results in battery constraints with more complex non-linear terms that depend on both payload weight and arrival times).

To extend the formulation for this case, a new set of non-negative variables $r_v$ are defined to represent the time that a drone arrives at vertex $v$ in minutes; $v \in V$ (note that this can differ from the time that they serve the customer). Moreover, the non-negative parameter $eh_v$ is introduced, which equals the energy consumed in each unit of time while the drone is hovering/waiting at vertex $v \in V$ of the generated graph, which are determined in the preprocessing phase (recall that in the generated graph, in contrast to the original graph, the load of the drone is specified in each vertex, and thus the energy consumed per time unit when the drone hovers/waits can be determined in the preprocessing phase). Then, objective (31) and constraints (35) must be slightly changed to include the extra energy consumed while drones arrive sooner as follows:

$$\min \sum_{a \in A} c_a \times x_a + \delta \sum_{a \in A} ed_a x_a + \delta \sum_{v \in V} eh_v \big(y_{cl(v)} - r_v\big) \tag{50}$$

subject to: (32)–(34), (36)–(49), and

$$r_v \leq y_{cl(w)} + t_a + M^7_{wv}(1 - x_a) \qquad a = (w,v) \in A: w \neq s, v \neq e \tag{51}$$

$$f_{cl(v)} - ed_a - eh_v \times \big(y_{cl(v)} - r_v\big) \geq f_{cl(w)} - M^4_a(1 - x_a) \qquad a = (w,v) \in A. \tag{52}$$

$$r_v \geq 0 \qquad v \in V. \tag{53}$$

Note that (35) is replaced by (52). Moreover, (51) is newly added to compute the exact arrival time of a drone at each vertex (notice the inequality sign in (51)). The objective (31) is also replaced by (50) to capture the additional costs of early presence of drones at customer locations. There are other approaches to formulate the extension, but they need to linearize many products of binary and continuous variables using big-M constraints.

### 4.4 Extension for covering load-dependent flight times and flight policies

A flight policy (determining all the controllable factors of the flight such as velocity and angle at each flight phase) is commonly considered independent of payload weight in the drone routing literature. In this case, the flight time remains constant and only the energy consumption varies by the payload weight. However, there may be a flight policy that depends on payload weight. For example, for heavier payloads, lower velocities are set in order to save more energy, have a safer flight, or consider the drones' technical limitations. In this case, both times and energy consumption for flying an arc can be dependent on the payload weight. Fortunately, the new formulation can be extended to cover this case.

Let us define the variable $z_{vw}$ (instead of $z_{ij}$) as a binary variable that is equal to 1 if the same drone performs the trip whose first request is in vertex $w$ immediately after a trip whose last request is in vertex $v$, and 0 otherwise, $v = (r, L) \in V^-: |L| = 1, w \in V^-: cl(w) \neq cl(v)$. Then, constraints (38), (39), (42), (43), and (44) must be replaced by the following constraint sets, respectively:

$$y_i + t_a \leq y_j + M_{ij}^8(1 - x_a) \qquad a \in A \tag{54}$$

$$y_{cl(v)} + t_{a=(v,e)} + t_{b=(s,w)} \leq y_{cl(w)} + M_{vw}^9(1 - z_{vw})$$
$$v = (r, L) \in V^-: |L| = 1, w \in V^-: cl(w) \neq cl(v) \tag{55}$$

$$\sum_{v=(r,L) \in V^-: |L|=1:\ cl(v) \neq cl(w)} z_{vw} \leq x_{a=(s,w)} \qquad w \in V^- \tag{56}$$

$$\sum_{w:\ cl(v) \neq cl(w)} z_{vw} \leq x_{a=(v,e)} \qquad v = (r, L) \in V^-: |L| = 1 \tag{57}$$

$$\sum_{a \in \delta^{out}(s)} x_a - \sum_{v=(r,L) \in V^-: |L|=1} \sum_{w \in V^-:\ cl(w) \neq cl(v)} z_{vw} \leq N \tag{58}$$

The extended optimization model again has a linear structure because the time and energy required for flying each arc with each set of requests can be determined in the preprocessing phase for load-dependent flight polices.

# 5. Computational Experiments

This section presents our numerical study to evaluate the computational performance of the new formulation of the DRP proposed in Section 4.2. Section 5.1 introduces the benchmark instances (test problems) and Section 5.2 compares the new formulation and the most recent algorithm proposed for solving the old formulation of the DRP with the special energy functions (25) and (26). The details of the computational study are given in Appendix D.

## 5.1. Benchmark instances

The performance of the proposed formulation procedure is evaluated on benchmark instances recently used by Cheng et al. (2020), which are generated based on older instance-generation frameworks used in Solomon (1987) and Dorling et al. (2017). For type-1 instances, included in set A1, the depots are located at the lower left corner of the district. For type-2 instances, included in set A2, the depots are located in the middle of the area. The number of the requests in these instances ranges from 10 to 50. Cheng et al. (2020) reported the results considering three different sets of the parameters appearing in the objective function. They used symbols R, E, and R + E to represent three objectives that minimizes only the travel cost ($\delta = 0$), only the energy cost ($c_a = 0$, for all $a \in A$), and both travel and energy costs, respectively (these are here called *objective settings*). All the parameters and optimization settings are set similarly as in Cheng et al. (2020) (e.g. the experiments are conducted on the same number of threads specified in their paper).

Table 1 summarizes the metadata of the benchmark instances. The first column in this table identifies the number of requests for each instance. The next columns present the value of parameter $UPMNR$ (used for the pre-analysis of the generated graph's dimensions, see Appendix B), the number of vertexes, and the number of arcs of the generated graph (Section 4.1.1), respectively. The last two columns specify upper bounds (derived in the pre-analysis described in Section 4.1.4) on the cardinalities of the vertex and arc sets of the generated graph using the equations (29) and (30), respectively. By comparing the actual and worst-case number of vertexes/arcs (e.g., compare 422,214 with 1,536,315,987 for instance 1-45-1), one can see that the generated graphs in practice are much smaller than what expected in the worst case.

In Figure 4, the boxplots indicate spreads of the ratios of the number of vertexes and the number of arcs of the generated graph to their upper bounds. This figure also shows the ratios of the time spent in the preprocessing stage to the total execution time (including solver time) for the three objective settings. In all the boxplots, the median is closer to the bottom of the box and the whisker is shorter on the lower end of the box, meaning all the ratio distributions are positively skewed. The figure confirms that the third quartiles and data concentrations are less than 5% for all the ratios. Hence, the size of the generated graph is much smaller than what expected in the worst case. Moreover, the times required for producing the generated graph in the preprocessing phase are considerably small compared to the times required by the optimization solver.

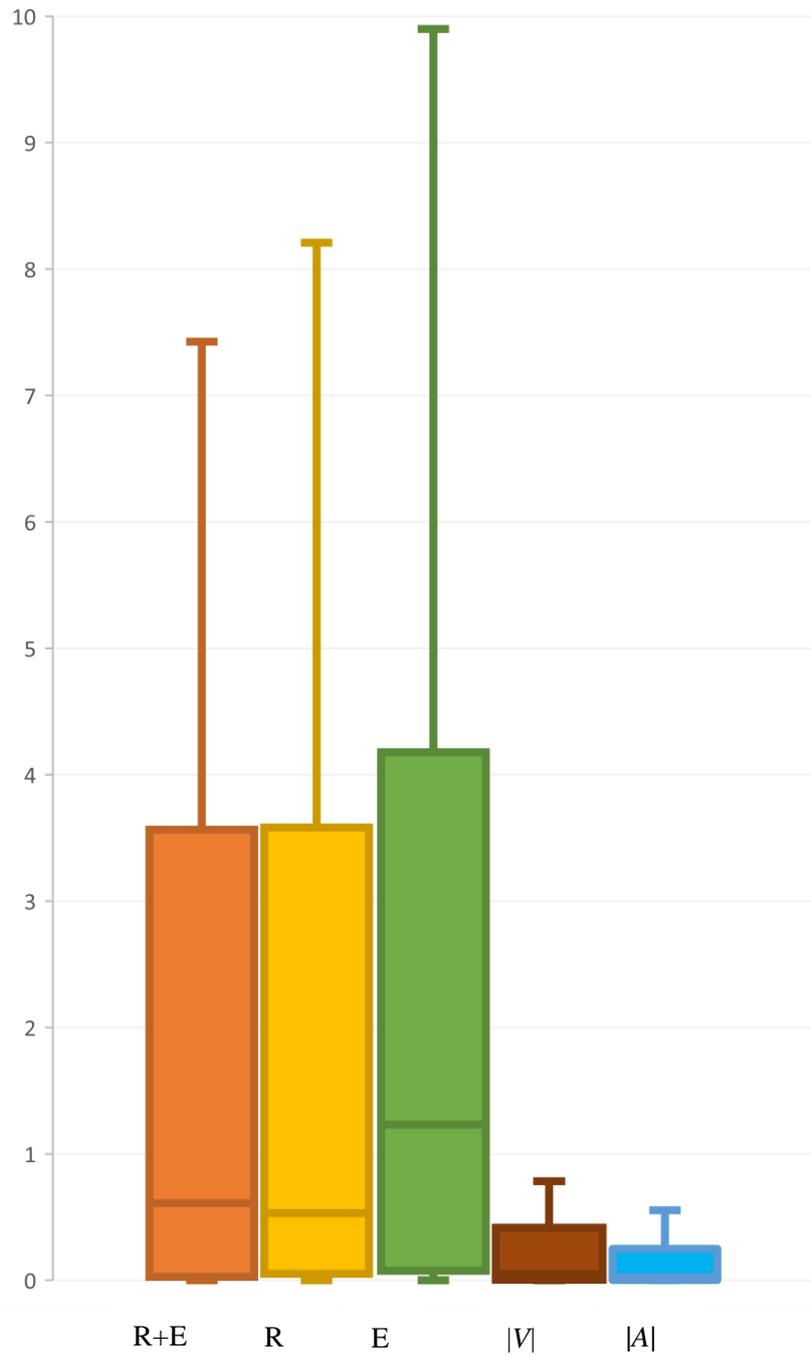

**Figure 4.** Spreads of the ratios of the number of vertexes and the number of arcs of the generated graph to their corresponding upper bounds given in equations (29) and (30), and the ratios of the time spent in the preprocessing stage to the total execution time (including solver time) for the three objective settings (all numbers are percentages).

**Table 1.** Metadata of benchmark instances: name, number of customers, UPNMR used in (29) and (30) in pre-analysis of generated graph's size, generated graph's actual dimensions, and upper bounds on generated graph's dimensions

| Instance | |R| | UP-MNR | |V| | |A| | Upper bound |V| | Upper bound |A| | Instance | |R| | UP-MNR | |V| | |A| | Upper bound |V| | Upper bound |A| |
|---|---|---|---|---|---|---|---|---|---|---|---|---|---|
| 1-10-1 | 10 | 4 | 156 | 400 | 1302 | 4760 | 2-10-4 | 10 | 4 | 101 | 223 | 1302 | 4760 |
| 1-10-2 | 10 | 5 | 148 | 430 | 2562 | 11186 | 2-10-5 | 10 | 5 | 83 | 187 | 2562 | 11186 |
| 1-10-3 | 10 | 4 | 204 | 588 | 1302 | 4760 | 2-15-1 | 15 | 7 | 397 | 1103 | 97142 | 603416 |
| 1-10-4 | 10 | 5 | 334 | 1034 | 2562 | 11186 | 2-15-2 | 15 | 5 | 241 | 571 | 22067 | 102916 |
| 1-10-5 | 10 | 5 | 170 | 492 | 2562 | 11186 | 2-15-3 | 15 | 5 | 173 | 405 | 22067 | 102916 |
| 1-15-1 | 15 | 5 | 674 | 2088 | 22067 | 102916 | 2-15-4 | 15 | 7 | 291 | 731 | 97142 | 603416 |
| 1-15-2 | 15 | 8 | 6659 | 29347 | 148622 | 1018688 | 2-15-5 | 15 | 4 | 194 | 468 | 7052 | 26840 |
| 1-15-3 | 15 | 6 | 2074 | 7582 | 52097 | 285098 | 2-20-1 | 20 | 6 | 329 | 801 | 333282 | 1888184 |
| 1-15-4 | 15 | 5 | 918 | 2892 | 22067 | 102916 | 2-20-2 | 20 | 6 | 854 | 2460 | 333282 | 1888184 |
| 1-15-5 | 15 | 6 | 952 | 3116 | 52097 | 285098 | 2-20-3 | 20 | 6 | 614 | 1708 | 333282 | 1888184 |
| 1-20-1 | 20 | 6 | 1385 | 4697 | 333282 | 1888184 | 2-20-4 | 20 | 4 | 405 | 967 | 23202 | 89720 |
| 1-20-2 | 20 | 6 | 1315 | 4121 | 333282 | 1888184 | 2-20-5 | 20 | 7 | 2455 | 8743 | 875922 | 5713796 |
| 1-20-3 | 20 | 5 | 1845 | 6189 | 100722 | 481196 | 2-25-1 | 25 | 6 | 1059 | 2885 | 1386377 | 7983630 |
| 1-20-4 | 20 | 7 | 4861 | 19552 | 875922 | 5713796 | 2-25-2 | 25 | 7 | 1381 | 4131 | 4751277 | 31672526 |
| 1-20-5 | 20 | 7 | 2160 | 7093 | 875922 | 5713796 | 2-25-3 | 25 | 6 | 1117 | 3137 | 1386377 | 7983630 |
| 1-25-1 | 25 | 6 | 3804 | 12658 | 1386377 | 7983630 | 2-25-4 | 25 | 6 | 855 | 2261 | 1386377 | 7983630 |
| 1-25-2 | 25 | 9 | 15186 | 63901 | 31790652 | 267455876 | 2-25-5 | 25 | 6 | 806 | 2092 | 1386377 | 7983630 |
| 1-25-3 | 25 | 6 | 4985 | 17007 | 1386377 | 7983630 | 2-30-1 | 30 | 6 | 1409 | 3869 | 4397882 | 25561936 |
| 1-25-4 | 25 | 7 | 3341 | 11483 | 4751277 | 31672526 | 2-30-2 | 30 | 7 | 3013 | 9195 | 18648482 | 125791156 |
| 1-25-5 | 25 | 8 | 6211 | 24083 | 13403877 | 101239430 | 2-30-3 | 30 | 7 | 4572 | 15875 | 18648482 | 125791156 |
| 1-30-1 | 30 | 7 | 27714 | 120222 | 18648482 | 125791156 | 2-30-4 | 30 | 8 | 3738 | 11836 | 65471882 | 501939136 |
| 1-30-2 | 30 | 6 | 14167 | 56101 | 4397882 | 25561936 | 2-30-5 | 30 | 8 | 4886 | 17028 | 65471882 | 501939136 |
| 1-30-3 | 30 | 7 | 10982 | 40436 | 18648482 | 125791156 | 2-35-1 | 35 | 8 | 8521 | 30638 | 246950622 | 1910298332 |
| 1-30-4 | 30 | 7 | 19836 | 82951 | 18648482 | 125791156 | 2-35-2 | 35 | 8 | 1526 | 4056 | 246950622 | 1910298332 |
| 1-30-5 | 30 | 8 | 23211 | 100509 | 65471882 | 501939136 | 2-35-3 | 35 | 10 | 7106 | 25046 | 2718211722 | 26058049652 |
| 1-35-1 | 35 | 7 | 80877 | 363253 | 58664062 | 398626236 | 2-35-4 | 35 | 9 | 3481 | 11257 | 882417762 | 7647658796 |
| 1-35-2 | 35 | 9 | 50879 | 227720 | 882417762 | 7647658796 | 2-35-5 | 35 | 9 | 7965 | 28365 | 882417762 | 7647658796 |
| 1-35-3 | 35 | 9 | 37475 | 157790 | 882417762 | 7647658796 | 2-40-1 | 40 | 8 | 6939 | 22677 | 772459522 | 6010931968 |
| 1-35-4 | 35 | 7 | 21821 | 83593 | 58664062 | 398626236 | 2-40-2 | 40 | 9 | 5906 | 19162 | 3233409442 | 28221004996 |
| 1-35-5 | 35 | 8 | 23085 | 97586 | 246950622 | 1910298332 | 2-40-3 | 40 | 8 | 5251 | 16941 | 772459522 | 6010931968 |
| 1-40-1 | 40 | 7 | 30647 | 122297 | 157222042 | 1073651191 | 2-40-4 | 40 | 8 | 4989 | 15905 | 772459522 | 6010931968 |
| 1-40-2 | 40 | 7 | 28943 | 118464 | 157222042 | 1073651191 | 2-40-5 | 40 | 9 | 13725 | 47431 | 3233409442 | 28221004996 |
| 1-40-3 | 40 | 8 | 63337 | 275501 | 772459522 | 6010931968 | 2-45-1 | 45 | 8 | 10711 | 38007 | 2097702632 | 16391903244 |
| 1-40-4 | 40 | 8 | 34164 | 139298 | 772459522 | 6010931968 | 2-45-2 | 45 | 7 | 21763 | 87842 | 373277072 | 2558178196 |
| 1-40-5 | 40 | 7 | 15286 | 55986 | 157222042 | 1073651191 | 2-45-3 | 45 | 6 | 6754 | 21282 | 55619732 | 327517764 |
| 1-45-1 | 45 | 11 | 422214 | 2194943 | 153631598717 | 1638779370802 | 2-45-4 | 45 | 8 | 21498 | 79276 | 2097702632 | 16391903244 |
| 1-45-2 | 45 | 8 | 68913 | 304543 | 2097702632 | 16391903244 | 2-45-5 | 45 | 8 | 3238 | 9426 | 2097702632 | 16391903244 |
| 1-45-3 | 45 | 8 | 76252 | 336828 | 2097702632 | 16391903244 | 2-50-1 | 50 | 8 | 23908 | 92662 | 5101140502 | 39985492360 |
| 1-45-4 | 45 | 8 | 68163 | 287747 | 2097702632 | 16391903244 | 2-50-2 | 50 | 9 | 10973 | 37231 | 27650043802 | 243376600126 |
| 1-45-5 | 45 | 10 | 65080 | 292364 | 41975043707 | 408076008914 | 2-50-3 | 50 | 8 | 13226 | 42175 | 5101140502 | 39985492360 |
| 2-10-1 | 10 | 5 | 127 | 317 | 2562 | 11186 | 2-50-4 | 50 | 9 | 12548 | 42990 | 27650043802 | 2.43377E+11 |
| 2-10-2 | 10 | 4 | 73 | 161 | 1302 | 4760 | 2-50-5 | 50 | 8 | 7821 | 25505 | 5101140502 | 39985492360 |
| 2-10-3 | 10 | 5 | 133 | 323 | 2562 | 11186 | | | | | | | |

## 5.2. Computational performance of new formulation approach on existing test problems

The algorithm used in the preprocessing phase, the new formulation, and Cheng et al.'s branch-and-cut algorithm were all coded in C++, and the MILP formulations were solved by CPLEX 12.8. Computations were performed on a computer with an Intel(R) Core(TM) i5-7400 and 3 GHz CPU and 16GB of RAM, running on a 64-bit Windows operating system. All the parameters are set to their default values in the solver except specified. All times are in seconds.

Cheng et al. (2020) applied a branch-and-cut algorithm, called *CAR algorithm* hereafter, to solve the old formulation (Section 3.3) efficiently on a cluster of Intel Xeon X5650 CPUs with 2.67 GHz and 24 GB RAM under Linux 6.3. Here, we report the results by directly solving the new formulation with CPLEX while using the convex energy-consumption formula used in Cheng et al. (2020). In some tables, Cheng et al. (2020) initiate the solving procedure with the solution found in other models to accelerate the search procedure, whereas in our setting no appropriate initial solution is used. Cheng et al. (2020) used a single thread processor setting in general solver and a run time limit of 4 hours for small-size instances ($|R| <= 30$) and multi-thread processor setting and 12 hours run-time limit for the ones with larger sizes ($|R| > 30$). Here, the same number of threads was enabled in the solver. In results reported by Cheng et al. (2020), some instances are solved through both settings. Here, we consider the setting with the better quality in their paper. Details of the solutions and performance of our formulation for the instances regarding the three objective values are reported and discussed in Appendix D.

The stopping criterion for Cplex is set to 0.01% relative error bound, as in Cheng et al. (2020). Since the CUP of our computer is 3 GHz, while it is 2.67 GHz in Cheng et al. (2020), all of our running times are multiplied by 3/2.67 to have a fair comparison. The run times and bounds reported in Cheng et al. (2020) are almost similar to the ones obtained based on our implementation of the Cheng et al.'s algorithm, but their results are better in some cases, so their results are used here for comparison in favor of their algorithm (we could not run their code as it was not made available online). However, note that the differences in run times of the old and new methods are extremely large (the former run times are 20 times larger than the latter on average), and thus slight changes in the run times cannot affect the result of our comparison.

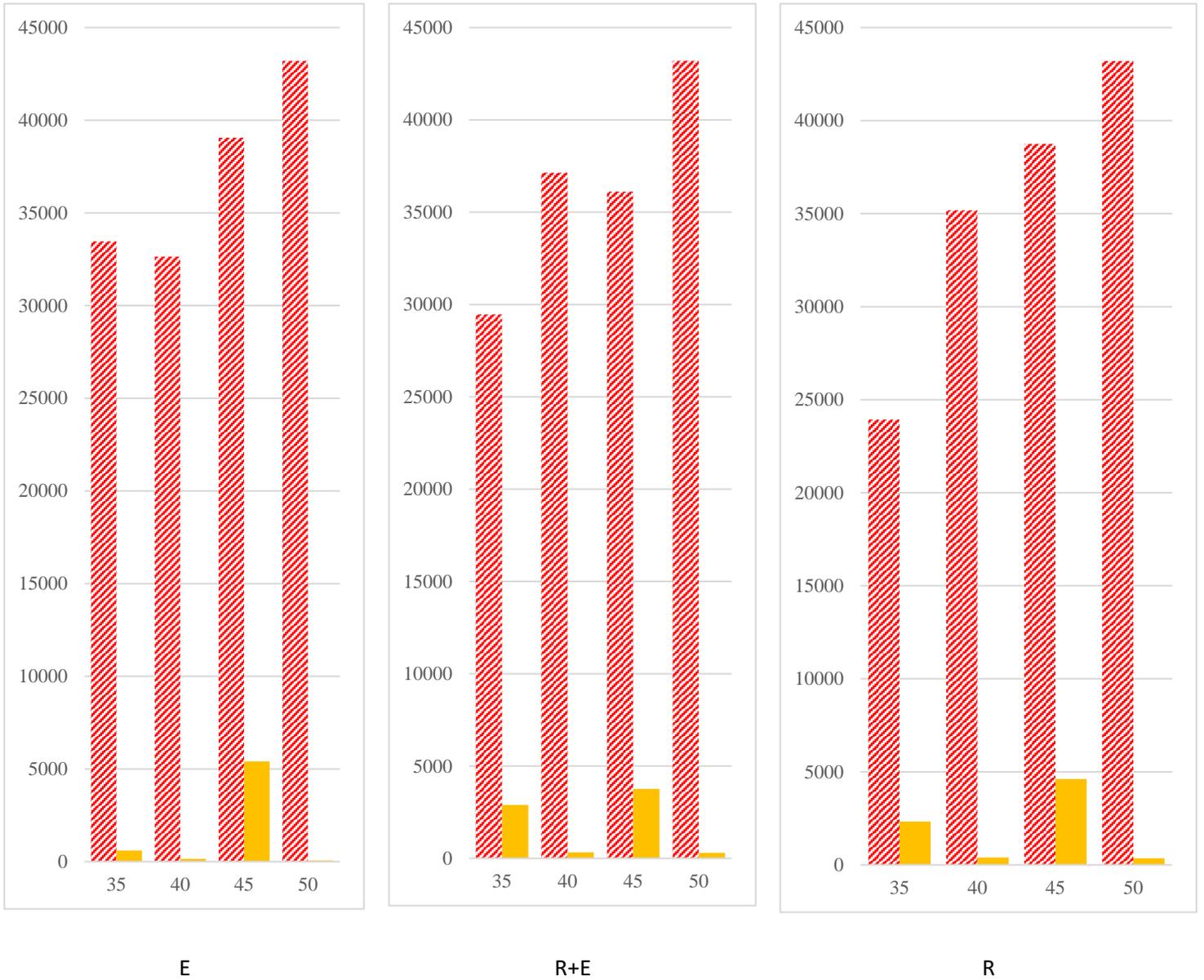

**Figure 5.** Average run times (including preprocessing and solver run times) for solving instances with 35, 40, 45, and 50 customers for three objective settings; ▨ : CAR algorithm, ▮ New formulation

Figure 5 compares average the run times (including preprocessing and solver run times) consumed for solving large-size instances using the new formulation with the those reported by Cheng et al. (2020). Figure 6 illustrates the ratio of the total run time (including preprocessing and solver run time) consumed for solving each instance to the reported total run time by Cheng et al. (2020). Recall that the experiments are conducted on the same number of threads specified in Cheng et al. (2020). From these figures one can see the following observations.

- In Cheng et al. (2020), for 71 out of 255 instances, the 0.01%-optimality of the obtained solution is not proved, while *all* instances are solved to 0.01%-optimality using the new formulation (our experiments indicate that all instances can be solved optimally by the new formulation within the time limit).

- The results approve that the total run times consumed for solving the instances are reduced significantly using the new formulation. The average run times of the CAR algorithm for solving the objective settings R+E, R, and E are 15752, 19065, and 16541 seconds, respectively; whereas the similar run times (including both preprocessing and solver times) of the new formulation drastically reduce to 835, 1144 and 706 seconds.

Our results indicate that the run time exceeds 3600 seconds (one hour) only for 8 instances, while this happens for 130 instances in Cheng et al. (2020).

An explanation may be needed about the non-increasing order of the average run times consumed by the new formulation for solving instances with 35 and 40, 45, and 50 customers, depicted in Figure 5, which is similarly observed for the three objective settings R, E, and R + E. According to Table 1, the average cardinalities of the vertex sets of the generated graphs for instances with 35, 40, 45 and 50 customers are 24,273.6, 20,918.7, 76,458.6 and 13,695.2, which can be ranked as 2nd, 1st, 3rd, and 4th, respectively (recall that generated graphs are the same for all the three objective settings R, E, and R + E). Why this happens is out of our control as all the instances are used from the literature of DRP without any modifications (it seems that the way used for generating the instances could not control the computational difficulty of the instances as the size increases; indeed, an instance of routing problem with smaller size is not necessarily easier than one with a larger size). Furthermore, note that for instances of the set A1, the number of vertexes and arcs of the generated graphs are much greater in compare to instances of set A2 with the same number of customers. This happens because when the service area is larger, many of the potential vertexes are eliminated in the preprocessing stage as many points are as far from each other that a drone cannot fly safely between the points. Considering the non-increasing order of the average sizes of the generated graph (given above), one can see that the average run times reported in Figure 5 are sensible because the average run time increases as the average size rises.

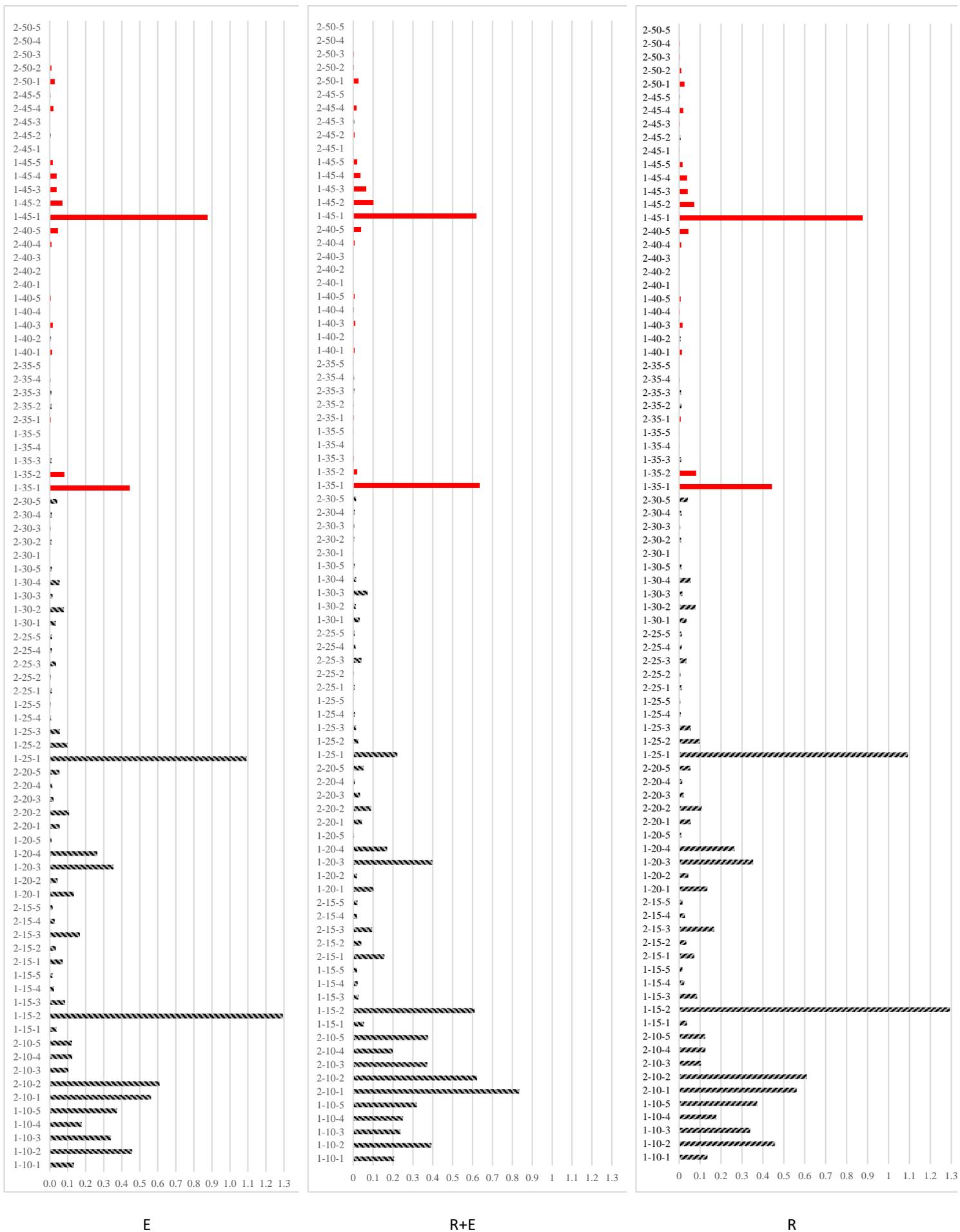

**Figure 6.** Ratio of the total run times (including preprocessing and solver run times) of the new formulation to the run time of the CAR algorithm. Red bars represent those instances for which the CAR algorithm cannot find a 0.01%-optimal solution within the run time limit (note that for all the instances the new formulation finds desired solutions quickly, never passing the run time limit).

For the objective settings R+E and R, Cheng et al. (2020) also compare the run times of solving the 50 small-size instances ($|R| \leq 30$) with the run times when the linear energy-consumption formula is used instead of the convex function (which can be handled more easily by their algorithm). As expected, for the linear case, the run times of the CAR algorithm are much smaller. However, according to Figure 7, they are still very large compared to the run times of our new formulation.

In sum, our numerical results clearly show that solving the new formulation is computationally efficient, while it has the capability to adopt any energy-consumption model.

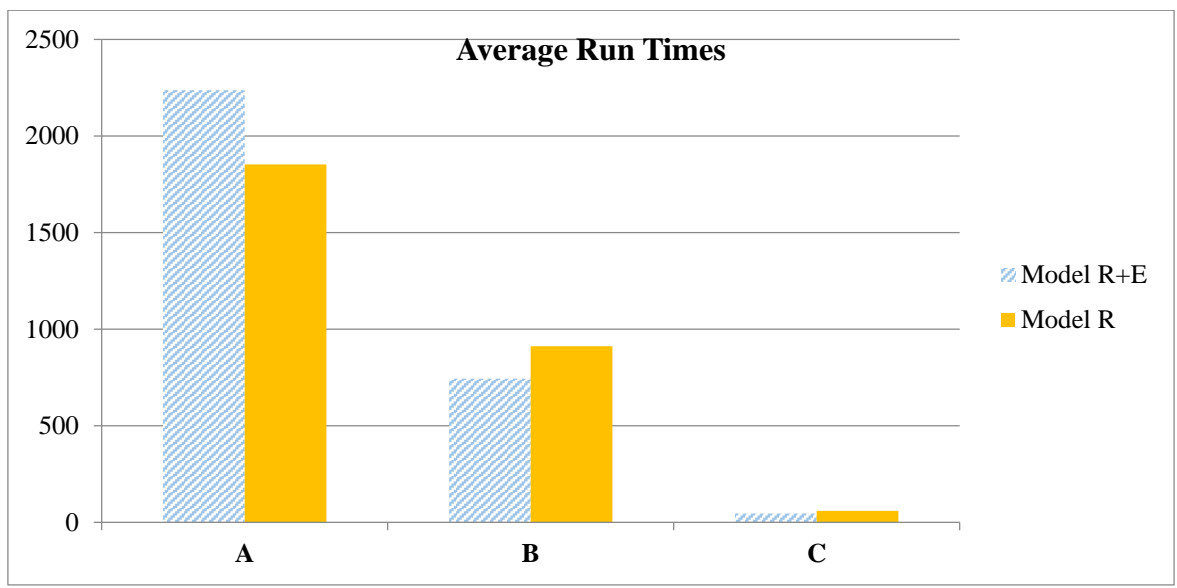

**Figure 7.** Average run times for solving small-size instances. **A:** CAR algorithm + Convex energy-consumption formula, (25); **B:** CAR algorithm + Linear energy-consumption formula, (26); **C:** New formulation + Convex energy consumption formula, (25)

### 5.3 Testing on new instances with up to 150 customers

With the purpose of exploring the performance of the new formulation on larger-size test problems, instances with a larger number of customers are generated in the same way described in Cheng et al (2020). The new generated test problems have 100 to 150 customers scattered in the service area and the results obtained using the new solution procedure for the objective setting R+E are presented in Table 2. In this table, the columns UP, LB, 4 Hours Gap%, 2 Hours Gap%, Solver Time and preprocessing Time, represent the upper bound on the objective value reported by the solver, the lower bound on the objective value reported by the solver, relative error bound after 4 hours of run time reported by the solver, relative error bound after 2 hours of run time reported by the solver, the time spent in solver and the time spent in the preprocessing stage to construct the generated graph for the instance, respectively.

A run time limit of 4 hours is set in the solver's setting and for better performance of the solution process, the following constraints and valid inequalities are added to model.

$$y_j \geq \sum_{a=(v,w)\in A: cl(v)=i, cl(w)=j} x_a t_{ij} \quad j \in R \tag{59}$$

$$y_j \geq \sum_i z_{ij}(t_{i,n+1} + t_{0j}) \quad j \in R \tag{60}$$

$$f_j \geq \sum_{a \in \delta^{in}(v), v \in V_j} ed_a x_a \quad j \in R \tag{61}$$

$$f_j \leq M - \min_{v: a=(v,e), v \in V_j} ed_a \quad j \in R \tag{62}$$

$$u_j \geq u_i + 1 - |R|\left(1 - \sum_{a=(v,w)\in A: cl(v)=i, cl(w)=j} x_a\right) \quad i, j \in R \tag{63}$$

$$nu = \sum_{a \in \delta^{out}(s)} x_a \tag{64}$$

Lower bounds of the time at which a drone starts providing service to a customer is determined through constraints (H.1) and (H.2). These constraints state that the service time of the customer is greater than the time required for flying the last arc and the last two arcs on the path for reaching the customer, respectively. Furthermore, lower and upper bounds for the accumulated energy consumption of a drone upon arrival at a customer's place is determined using valid inequalities (H.3) and (H.4). Constraints (H.3) are derived similarly to Constraints (H.1) where the energy is considered instated of time. Constraints (H.4) are valid because the drone must have at least the minimum energy needed to fly back to the depot when visiting any customer to prevent drone crash. In this constraint set, the parameter $\min_{v: a=(v,e), v \in V_j} ed_a$ equals the minimum energy needed to fly from customer $j$ to the end depot if flights energy consumption satisfy triangle inequality (if not, it must be replaced with the result obtained through solving a shortest path problem, see Appendix B for more details). Constraints (H.5) are adaption of MTZ sub-tour elimination constraints introduced in Miller, Tucker, and Zemlin (1960). In these constraints, the non-negative variables $u_i$, $i \in R$, signify the order of the customers visited in each trip. Furthermore, an additional auxiliary *integer* variable named $nu$ is set to the number of trips completed by the drones in constraint (H.6), which enables the MILP solver branches on the values of $nu$ to increase its speed.

Table 2. The results obtained through solving large size instances with setting R+E using new formulation.

| Instance | UP | LB | 4 Hours Gap% | 2 Hours Gap% | Solver Time | Preprocessing Time |
|---|---|---|---|---|---|---|
| A1_100_1 | 41595.08 | 41595.08 | 0.00 | 0.00 | 7023.16 | 923.25 |
| A1_100_2 | 42848.06 | 41713.75 | 2.67 | 2.81 | 14400 | 749.75 |
| A2_100_1 | 44008.96 | 43953.13 | 0.13 | 0.20 | 14400 | 58.36 |
| A2_100_2 | 47440.53 | 47418.19 | 0.05 | 0.07 | 14400 | 17.21 |
| A2_150_1 | 66132.33 | 66065.38 | 0.10 | 0.28 | 14400 | 372.09 |
| A2_150_2 | 68527.18 | 66968.67 | 2.27 | 2.40 | 14400 | 622.05 |

As Table 2 indicates, the maximum error bound for any instance is around 3%. The results confirm that the larger number of customers increases the total cost and energy consumption as expected. The same as the results obtained from solving the benchmark test problems, the CPU times reported are much less for set A2 instances where the service area is larger and the depot is located at the middle of the area. In terms of computational performance, our algorithms could generally solve larger size instances compared to those considered in exact algorithms for the DDP despite the fact that our formulation can use a wide range of nonlinear and more precise energy consumption formulas without adding to the complexity of the problem.

## 6. Concluding remarks

Drones have special characteristics that make a delivery problem with drones different from a similar problem in which ground vehicles are used. Drones have major restrictions on their maximum flight durations and payload capacities. Seeing that these restrictions considerably shrink the solution space of a drone delivery problem, this paper proposes a formulation approach that constructs a linear optimization model whose structure is independent of how the drones' energy consumption is expressed. The new approach enables us to use the best available energy consumption model, while the existing studies either assume that flight range is a fixed value or use some approximations (that are linear or convex in the weight of the drone and its payload) of the actual energy-consumption models. In fact, the new approach computes the energy consumed for carrying any reasonable set of requests in each drone fly between two specific points as an input, which can be determined through any computational method, e.g., using explicit formulas, simulation, or black-box procedures. Hence, the effects of the drone's and parcel's weight and size, the wind direction and speed, the battery weight, the weather conditions, the value and angle of drone's speed and acceleration, the weight and size of the payload on the energy consumption, and the energy needed for different drone's motions (i.e., takeoff, landing, hovering, and horizontal flight) in delivery trips and empty returns can be considered.

An extensive numerical study is carried out to evaluate the computational performance of the new formulation on 255 benchmark instances, which were recently solved by a B&C algorithm. For 71 out of the 255 instances, the B&C algorithm could not find any 0.01%-optimal solutions within the preset time limit, while all the 255 instances are solved to 0.01%-optimality using the new formulation (all instances are also solved exactly within the time limit). The total run times required for using the new formulation are reasonably smaller when compared to the B&C algorithm. The average run times required by the B&C algorithm are 15752, 19065 and 16541 seconds for the objective settings R+E, R and E, respectively; while they reduce to 835, 1144, and 706 seconds, respectively. In addition, while instances in the relevant literature have at most 50 customers, new instances with up to 150 customers are generated and solved using the new solution procedure. This shows that the new formulation not only has the capability of incorporating arbitrary energy-consumption formulas, it is also computationally efficient.

The formulation procedure presented in this paper may be effectively applied to other drone routing problems. Here, it is shown how it can be used for linearly modeling two challenging cases where the drones cannot land out of the time windows and where both energy consummation rate and time of a fly depend on the payload weight. In future studies, the new approach can be used for addressing more complex settings, for example, where a company uses a fleet of heterogeneous drones or a hybrid fleet of drones and trucks.

Stochastic or time-related variations can also impact a drone's energy consumption. However, as the flight range of a drone is typically very limited, one can assume such variations are small and can be handled by a reserve battery capacity. It may be the main reason that most of the papers on drone routing are established in static and deterministic settings. Another reason can be the difficulty of solving the existing drone routing problems (before this study no study can solve benchmark instances of the DRP with 50 customers exactly). Studying the DRP considered in this paper in a stochastic or time-dependent setting remains as a future study for relevant applications, in which such variations are very influential even in a short period, or the drones' ranges are long and such variations cannot be covered by the reserve battery capacity. In fact, the linear structure of the new formulation that also linearly depends on the input parameters that can be uncertain increases the chance of successful application of existing stochastic or rubout optimization methods developed for MILPs.

## Appendix A. Parameters of energy consumption model (27)

For energy–consumption model (27), the parameters $\alpha_j^i$, $j = 1, \ldots, 14$, can be represented based on the parameters originally used by Kirschstein (2020) as follows:

$$\alpha_1^i = 0.5\rho. v^3. A. c_{Air} + P^{int}$$

$$\alpha_2^i = \kappa$$

$$\alpha_3^i = 0.25\rho^2. v^4. A^2. c_{Air}^2 + \rho. v^2. A. c_{Air} P^{climb}$$

$$\alpha_4^i = g^2$$

$$\alpha_5^i = \frac{27 c_{bd}\sqrt{g^3}}{\sqrt{\bar{c} n_{rotor} n_{blade} \rho r (\bar{c}_l)^3}}$$

$$\alpha_6^i = \frac{v^2 c_{bd}\sqrt{6 g \rho r n_{rotor} n_{blade} \bar{c}}}{4\sqrt{\bar{c}_l}}$$

$$\alpha_7^i = g. v. \sin\gamma$$

$$\alpha_8^i = \frac{1}{2\rho r^2 \pi n_{rotor}}$$

$$\alpha_9^i = v^2$$

$$\alpha_{10}^i = \frac{-\rho. v^3. A. c_{Air}}{g. \cos\gamma}$$

$$\alpha_{11}^i = -2. v. \tan\gamma$$

$$\alpha_{12}^i = (1 + \tan\gamma^2)$$

$$\alpha_{13}^i = \left(\frac{0.5. \rho. v^2. A. c_{Air}}{g. \cos\gamma}\right)^2$$

$$\alpha_{14}^i = \rho. v^2. A. c_{Air}. \tan\gamma$$

## Appendix B. Complementary details on rules used for reducing generated graph's dimensions

To reduce the number of vertices of the generated graph, using the rules R1, R2, and R3 explained in Section 4.1, the payload capacity constrains are fully checked and the time-window and energy constraints are partially checked before creating each vertex (recall that this partial checking is optional, but very effective in downsizing the generated graph). The rule R1 always works properly. The rules R2 and R3 also work properly and do not eliminate any feasible solution if the following two conditions hold:

    A. The travel times are metric (i.e., the triangular inequality is satisfied for all of the travel times).

    B. The energy required for serving each set of customers cannot be less than the energy required for serving any *subset* of that set.

Though these conditions are reasonable (and hold for benchmark instances solved in this paper; see below for a proof), but in some (pathological) cases they may be violated. In this situation (lack of the above conditions), the rule R2 must be modified as follows to ensure that no feasible solution is removed:

R2'. Partial checking of the battery and time-window constraints: <u>for at least one permutation</u> of the requests contained in the set $L$, serving the requests in the order of the permutation is feasible in a single trip where the drone is not required to visit to any customer before arriving at the location of request $r$ from the depot with the following considerations:

i- in the partial checking of the time-window constraints: the traveling time between each pair of vertexes is replaced by the new time computed by solving a shortest path problem in the original graph whose arc weights are the original times (if condition A holds, the replacement of traveling times is not required);

ii- In the partial checking of the battery constraints: the energy required for traveling from a vertex for any given $S \subseteq L$ to the next vertex (determined in the permutation) is computed by solving a shortest path problem in the original graph where the weight of each arc is set to the energy required for flying the arc with the requests in $S$ (if condition B holds, the energy for flying directly to the next vertex with the requests in $S$ is only required).

$feas_{r,L} = 0$ if no permutation satisfying the above condition can be found.

It can be proven if condition (A) holds and the energy-consumption model is defined as the product of the travel time and energy-consummation rate per time unit (the energy functions (25) and (26) have such a structure), then condition B is also satisfied (as a result, both conditions A and B hold for the benchmark instances, in which the travel times are proportional to the Euclidian distances and energy functions (25) and (26) are used). To shows this, assume the energy-consumption model is given by

$$F_a(w_a, t_a) = t_a \theta(w_a) \qquad (B.1)$$

where all times are metric and $\theta(w_a)$ is an increasing function, denoting the energy required per unit of time for flying with payload weight $w_a$. Let us provide the proof for a small case that can be generalized to the general case. By contradiction, assume that condition B does not hold and the energy required for serving $L = \{1,2\}$ (visiting 1 after flying from the depot $e$) is strictly less than the energy required for serving $L = \{1\}$, that is,

$$E(e \to 1 \to 2 \to s) < E(e \to 1 \to s),$$

where the functional $E(\cdot)$ returns the energy for flying a given delivery path. Considering (B.1), this implies that

$$t_{e1}\theta(w_1 + w_2) + t_{12}\theta(w_2) + t_{2s}\theta(0) < t_{e1}\theta(w_1) + t_{1s}\theta(0)$$

Then, using $\theta(w_1) \leq \theta(w_1 + w_2)$ and $\theta(0) \leq \theta(w_2)$, one gets the following inequality:

$$t_{e1}\theta(w_1) + (t_{12} + t_{2s})\theta(0) < t_{e1}\theta(w_1) + t_{1s}\theta(0) \Leftrightarrow t_{12} + t_{2s} < t_{1s},$$

which contradicts the assumption that the travel times are metric ($t_{12} + t_{2s} \geq t_{1s}$). This completes the proof.

## Appendix C. An efficiently-computable upper bound on $MNR$

First recall that the parameter $MNR$ (or an upper on it) <u>is not used in any part of the new solution method</u>; it is <u>only used in a pre-analysis</u> of the size of the generated graph (before constructing it).

The exact value of the parameter $MNR$ is equal to the optimal value of the following mixed-integer non-linear programming model:

$$\max \sum_{a \in A} x_a - 1 \tag{C.1}$$

subject to:

$$\sum_{a \in \delta^{in}(v)} x_a \leq 1 \quad \forall v \in V_O^- \tag{C.2}$$

$$\sum_{a \in \delta^{out}(v)} x_a = \sum_{a \in \delta^{in}(v)} x_a \quad v \in V_O^- \tag{C.3}$$

$$\sum_{a \in \delta^{out}(s)} x_a = \sum_{a \in \delta^{in}(e)} x_a = 1 \tag{C.4}$$

$$y_v + t_a \leq y_w + M_a^2(1 - x_a) \quad a = (v, w) \in A_O \tag{C.5}$$

$$a_{r(v)} \leq y_v \leq b_{r(v)} \quad v \in V_O^- \tag{C.6}$$

$$a_d \leq y_v \leq b_d \quad v \in \{s, e\} \tag{C.7}$$

$$f_e \leq M \tag{C.8}$$

$$f_s = 0 \tag{C.9}$$

$$f_w - F_a(w_a, t_a) \geq f_v - M_a^1(1 - x_a) \quad a = (v, w) \in A_O \tag{C.10}$$

$$\sum_{a \in \delta^{out}(v)} w_a = \sum_{b \in \delta^{in}(v)} w_b - q_{r(v)} \quad v \in V_O^- \tag{C.11}$$

$$w_a \leq Q x_a \quad a \in A_O \tag{C.12}$$

$$\sum_{a \in \delta^{in}(e)} w_a = 0 \tag{C.13}$$

$$x_a \in \{0,1\} \quad a \in A_O \tag{C.14}$$

$$f_v \geq 0 \quad v \in V_O \tag{C.15}$$

$$y_v \geq 0 \quad v \in V_O \tag{C.16}$$

$$w_a \geq 0 \quad a \in A_O, \tag{C.17}$$

which cannot be solved simply. A naïve upper bound on $MNR$ is

$$\min\left\{|R|, \frac{Q}{\min_{r \in R}\{q_r\}}\right\},$$

which is very loose. An efficiently-computable tighter upper bound on $MNR$ for a given instance can be obtained by rounding down the optimal objective value of the following MILP model to an integer:

(C.1) subject to (C.2)–(C.9), (C.15), (C.16), and

$$f_w - F_a(q_{r(w)}, t_a) \geq f_v - M_a^1(1 - x_a) \qquad a = (v, w) \in A_O. \tag{C.18}$$

$$\sum_{a=(v,w) \in A} x_a \times q_{r(v)} \leq Q. \tag{C.19}$$

$$x_a \in [0,1] \qquad a \in A_O \tag{C.20}$$

This model is a relaxation of (C.1)– (C.17), in which variables $w_a$ and the constraints involving them are removed, constraints (C.10) are replaced by simple necessary conditions given in (C.18), and the payload capacity constraints (C.12) are replaced by (C.19). (C.20) is a relaxation of (C.14).

In Section 5.1, $UPMNR$ is set to this upper bound for computing the upper bounds on the number of vertexes and arcs of the generated graph in pre-analysis.

Note that one can solve the original model (C.2)–(C.9) by using our new formulation approach, but it is not reasonable here since the aim is to efficiently find a better upper on MNR instead of $\min\left\{|R|, \frac{Q}{\min_{r \in R}\{q_r\}}\right\}$ for comparing the generated graph's dimensions in theory and practice to assess the new formulation approach (this comparison implicitly shows how much the rules R1, R2, and R3 used for pruning the vertexes of the generated graph are effective).

## Appendix D. Details of numerical experiments in Section 5

This section provides the detailed results for the objective settings R+E, E ($c_a = 0$, for $a \in A$), and R ($\delta = 0$) in Table D.1, Table D.2, and Table D.3, respectively. Note that each table has 285 rows because some of the 255 instances were solved in more than one solver setting; in the summarizing figures, the best results obtained by the B&C algorithm are used for the sake of fairness.

In these tables, the second columns specify the single thread or multi thread setting in solver. Value 2 in the column Status indicates the instance is solved to optimality and value 9 states the solving procedure is stopped due to reaching the solver time limit. The next columns report the results of the B&C algorithm designed and implemented by Cheng et al. (2020). The columns UP, LB, Gap%, RLB, and Run Time represent the upper bound on the objective value reported by the solver, the lower bound on the objective value reported by the solver, relative error bound, the objective value in the solution of the relaxed problem and the total run time required for solving the instance, respectively. The next columns similarly describe the results obtained by solving the instances using the new formulation. The columns Solver Time and Preprocessing Time report the time spent in solver and the time spent in the preprocessing stage to construct the generated graph for each instance, respectively. When a better solution is found by the new method, the instance is highlighted with * in columns Status and UB.

**Table D.1.** The details results of experiments for the objective setting R+E.

| | | CAR algorithm | | | | | | New formulation | | | | | | |
|---|---|---|---|---|---|---|---|---|---|---|---|---|---|---|
| Instance | # of Threads | Status | UP | LB | Bound (%) | RLB | Run Time | Status | UP | LB | Bound (%) | RLB | Run Time | Solver Time | Preprocessing Time |
| 1-10-1 | 1 | 2 | 3133.64 | 3133.64 | 0.00 | 3011.79 | 0.88 | 2 | 3133.64 | 3133.64 | 0.00 | 2946.90 | 0.16 | 0.16 | 0.00 |
| 1-10-2 | 1 | 2 | 4740.72 | 4740.72 | 0.00 | 4695.62 | 0.19 | 2 | 4740.72 | 4740.72 | 0.00 | 4738.03 | 0.07 | 0.07 | 0.00 |
| 1-10-3 | 1 | 2 | 4557.82 | 4557.82 | 0.00 | 4555.37 | 0.45 | 2 | 4557.82 | 4557.82 | 0.00 | 4150.64 | 0.10 | 0.10 | 0.00 |
| 1-10-4 | 1 | 2 | 4392.71 | 4392.71 | 0.00 | 3886.16 | 1.28 | 2 | 4392.71 | 4392.71 | 0.00 | 3859.79 | 0.29 | 0.29 | 0.00 |
| 1-10-5 | 1 | 2 | 4526.00 | 4526.00 | 0.00 | 4349.29 | 0.40 | 2 | 4526.00 | 4526.00 | 0.00 | 4464.55 | 0.11 | 0.11 | 0.00 |
| 1-15-1 | 1 | 2 | 7074.75 | 7074.24 | 0.01 | 6741.96 | 12.79 | 2 | 7074.75 | 7074.75 | 0.00 | 6889.18 | 0.62 | 0.60 | 0.02 |
| 1-15-2 | 1 | 2 | 4399.81 | 4399.81 | 0.00 | 4051.32 | 22.51 | 2 | 4399.81 | 4399.81 | 0.00 | 4071.91 | 12.24 | 11.50 | 0.74 |
| 1-15-3 | 1 | 2 | 5970.07 | 5969.48 | 0.01 | 5165.03 | 1391.45 | 2 | 5970.07 | 5970.07 | 0.00 | 5334.69 | 33.34 | 33.25 | 0.09 |
| 1-15-4 | 1 | 2 | 5493.11 | 5492.81 | 0.01 | 5124.79 | 40.40 | 2 | 5493.10 | 5493.10 | 0.00 | 5360.05 | 0.81 | 0.81 | 0.00 |
| 1-15-5 | 1 | 2 | 7385.93 | 7385.48 | 0.01 | 6655.92 | 63.76 | 2 | 7385.93 | 7385.93 | 0.00 | 6759.71 | 1.15 | 1.14 | 0.01 |
| 1-20-1 | 1 | 2 | 8287.67 | 8287.67 | 0.00 | 8247.59 | 2.50 | 2 | 8287.66 | 8287.66 | 0.00 | 8204.83 | 0.23 | 0.22 | 0.01 |
| 1-20-2 | 1 | 2 | 9550.88 | 9550.79 | 0.00 | 9105.80 | 75.69 | 2 | 9550.88 | 9550.88 | 0.00 | 9368.38 | 1.34 | 1.33 | 0.01 |
| 1-20-3 | 1 | 2 | 8819.25 | 8819.25 | 0.00 | 8636.70 | 24.71 | 2 | 8819.25 | 8819.25 | 0.00 | 8575.07 | 8.77 | 8.74 | 0.03 |
| 1-20-4 | 1 | 2 | 6697.05 | 6696.54 | 0.01 | 6455.06 | 33.99 | 2 | 6697.05 | 6697.05 | 0.00 | 6590.14 | 5.19 | 4.94 | 0.25 |
| 1-20-5 | 1 | 2 | 7785.28 | 7784.52 | 0.01 | 7207.77 | 983.26 | 2 | 7785.28 | 7784.54 | 0.01 | 7419.15 | 3.85 | 3.80 | 0.05 |
| 1-25-1 | 1 | 2 | 10683.52 | 10682.45 | 0.01 | 10436.48 | 216.95 | 2 | 10683.52 | 10682.71 | 0.01 | 10188.98 | 42.86 | 42.75 | 0.11 |
| 1-25-2 | 1 | 2 | 8639.75 | 8638.89 | 0.01 | 8105.84 | 8600.06 | 2 | 8639.75 | 8639.51 | 0.00 | 8385.05 | 202.76 | 200.45 | 2.31 |
| 1-25-3 | 1 | 2 | 10097.80 | 10096.80 | 0.01 | 9774.80 | 1185.96 | 2 | 10097.80 | 10097.18 | 0.01 | 9663.37 | 15.88 | 15.68 | 0.21 |
| 1-25-4 | 1 | 2 | 10149.78 | 10148.81 | 0.01 | 9702.29 | 324.78 | 2 | 10149.78 | 10149.78 | 0.00 | 10007.64 | 2.67 | 2.57 | 0.10 |
| 1-25-5 | 1 | 9 | 11209.29 | 10899.57 | 2.76 | 10427.99 | 14400.01 | 2* | 11169.60* | 11169.40 | 0.00 | 10915.95 | 21.25 | 20.90 | 0.35 |
| 1-25-5 | 4 | 2 | 11169.40 | 11168.87 | 0.01 | 10427.99 | 8831.37 | 2 | 11170.09 | 11169.40 | 0.01 | 10915.95 | 9.76 | 9.41 | 0.35 |
| 1-30-1 | 1 | 2 | 9836.27 | 9835.30 | 0.01 | 9419.82 | 8584.35 | 2 | 9836.26 | 9836.26 | 0.00 | 9683.59 | 247.17 | 233.70 | 13.47 |
| 1-30-2 | 1 | 2 | 12669.33 | 12668.08 | 0.01 | 12361.69 | 1732.78 | 2 | 12669.33 | 12669.33 | 0.00 | 12599.82 | 21.19 | 19.67 | 1.53 |
| 1-30-3 | 1 | 9 | 12435.73 | 12077.86 | 2.88 | 11681.04 | 14400.00 | 2* | 12363.92* | 12363.39 | 0.00 | 11819.06 | 938.62 | 937.65 | 0.97 |
| 1-30-3 | 4 | 2 | 12363.92 | 12362.84 | 0.01 | 11681.04 | 23862.76 | 2 | 12363.92 | 12363.37 | 0.00 | 11819.06 | 304.00 | 285.95 | 0.97 |
| 1-30-4 | 1 | 2 | 12517.11 | 12515.88 | 0.01 | 12073.89 | 5211.59 | 2 | 12517.11 | 12516.18 | 0.01 | 12313.74 | 70.58 | 67.36 | 3.22 |
| 1-30-5 | 1 | 9 | 12090.40 | 12020.03 | 0.58 | 11644.46 | 14400.00 | 2 | 12090.40 | 12090.40 | 0.00 | 11932.47 | 34.25 | 29.45 | 4.80 |
| 1-30-5 | 4 | 2 | 12090.40 | 12089.71 | 0.01 | 11644.46 | 6318.78 | 2 | 12090.40 | 12090.40 | 0.00 | 11932.47 | 43.75 | 38.95 | 4.80 |
| 1-35-1 | 4 | 9 | 12510.16 | 12104.51 | 3.24 | 11693.53 | 43200.01 | 2* | 12439.34* | 12438.12 | 0.01 | 11921.93 | 24446.09 | 24374.15 | 71.94 |
| 1-35-2 | 4 | 9 | 13025.50 | 12663.33 | 2.78 | 12218.08 | 43200.01 | 2* | 13025.56 | 13024.35 | 0.01 | 12678.09 | 715.84 | 683.05 | 32.79 |
| 1-35-3 | 4 | 9 | 13236.10 | 13106.45 | 0.98 | 12706.73 | 43200.01 | 2* | 13236.10 | 13236.10 | 0.00 | 13207.88 | 182.44 | 171.00 | 11.44 |
| 1-35-4 | 4 | 9 | 13868.19 | 13838.51 | 0.21 | 13305.60 | 43200.01 | 2* | 13868.19 | 13868.19 | 0.00 | 13716.51 | 50.61 | 47.50 | 3.11 |
| 1-35-5 | 4 | 2 | 13286.73 | 13285.70 | 0.01 | 12709.65 | 20446.57 | 2 | 13286.73 | 13286.16 | 0.00 | 13126.90 | 35.11 | 31.35 | 3.76 |
| 1-40-1 | 4 | 9 | 15546.01 | 14931.91 | 3.95 | 14608.21 | 43200.01 | 2* | 15546.00 | 15546.00 | 0.00 | 15293.90 | 213.66 | 207.10 | 6.56 |
| 1-40-2 | 4 | 9 | 16887.24 | 16787.88 | 0.59 | 16156.12 | 43200.05 | 2* | 16887.21 | 16885.54 | 0.01 | 16741.74 | 77.82 | 72.20 | 5.62 |
| 1-40-3 | 4 | 9 | 14183.71 | 13648.06 | 3.78 | 13238.94 | 43200.01 | 2* | 14183.80 | 14182.44 | 0.01 | 14037.60 | 391.88 | 347.70 | 44.18 |
| 1-40-4 | 4 | 9 | 16292.14 | 16229.27 | 0.39 | 15562.99 | 43200.03 | 2* | 16292.14 | 16290.52 | 0.01 | 16213.09 | 134.43 | 124.45 | 9.98 |
| 1-40-5 | 4 | 9 | 15625.81 | 15264.11 | 2.32 | 14965.87 | 43200.03 | 2* | 15625.81 | 15624.26 | 0.01 | 15412.08 | 268.29 | 266.95 | 1.34 |
| 1-45-1 | 4 | 9 | 14760.22 | 14175.20 | 3.96 | 13817.26 | 43200.01 | 2* | 14575.02* | 14575.02 | 0.00 | 14284.45 | 23855.85 | 21471.90 | 2383.95 |
| 1-45-2 | 4 | 9 | 19747.19 | 19340.69 | 2.06 | 18942.58 | 43200.01 | 2* | 19734.46* | 19732.49 | 0.01 | 19321.99 | 3923.24 | 3871.25 | 51.99 |
| 1-45-3 | 4 | 9 | 18899.08 | 18434.22 | 2.46 | 17974.20 | 43200.01 | 2* | 18832.07* | 18831.50 | 0.00 | 18608.29 | 2504.56 | 2441.50 | 63.06 |
| 1-45-4 | 4 | 9 | 16341.14 | 15826.69 | 3.15 | 15398.76 | 43200.01 | 2* | 16304.98* | 16303.51 | 0.01 | 16018.90 | 1448.06 | 1395.55 | 52.51 |
| 1-45-5 | 4 | 9 | 18736.56 | 18307.23 | 2.29 | 17850.34 | 43200.01 | 2* | 18734.36* | 18733.31 | 0.01 | 18446.85 | 742.14 | 696.35 | 45.79 |
| 2-10-1 | 1 | 2 | 5000.23 | 5000.23 | 0.00 | 4772.59 | 0.23 | 2 | 5000.23 | 5000.23 | 0.00 | 4783.10 | 0.17 | 0.17 | 0.00 |
| 2-10-2 | 1 | 2 | 5827.05 | 5827.05 | 0.00 | 5825.59 | 0.12 | 2 | 5827.05 | 5827.05 | 0.00 | 5824.49 | 0.07 | 0.07 | 0.00 |
| 2-10-3 | 1 | 2 | 5271.53 | 5271.53 | 0.00 | 5267.79 | 0.20 | 2 | 5271.53 | 5271.53 | 0.00 | 5180.71 | 0.07 | 0.07 | 0.00 |
| 2-10-4 | 1 | 2 | 6158.59 | 6158.59 | 0.00 | 5401.12 | 0.80 | 2 | 6158.59 | 6158.59 | 0.00 | 5617.92 | 0.14 | 0.14 | 0.00 |
| 2-10-5 | 1 | 2 | 5535.34 | 5535.34 | 0.00 | 5112.11 | 0.34 | 2 | 5535.34 | 5534.94 | 0.01 | 5174.21 | 0.11 | 0.11 | 0.00 |
| 2-15-1 | 1 | 2 | 6871.81 | 6871.81 | 0.00 | 6868.38 | 0.95 | 2 | 6871.81 | 6871.81 | 0.00 | 6549.31 | 0.13 | 0.13 | 0.00 |
| 2-15-2 | 1 | 2 | 8537.57 | 8537.57 | 0.00 | 7542.77 | 10.08 | 2 | 8537.56 | 8537.56 | 0.00 | 7883.82 | 0.37 | 0.37 | 0.00 |
| 2-15-3 | 1 | 2 | 6614.52 | 6614.52 | 0.00 | 6465.54 | 1.46 | 2 | 6614.52 | 6614.52 | 0.00 | 6403.84 | 0.12 | 0.12 | 0.00 |
| 2-15-4 | 1 | 2 | 8780.18 | 8780.02 | 0.00 | 7878.04 | 9.95 | 2 | 8780.18 | 8780.18 | 0.00 | 7677.06 | 0.17 | 0.17 | 0.00 |
| 2-15-5 | 1 | 2 | 8674.76 | 8674.53 | 0.00 | 7978.28 | 5.83 | 2 | 8674.76 | 8674.76 | 0.00 | 8613.46 | 0.11 | 0.11 | 0.00 |
| 2-20-1 | 1 | 2 | 11425.87 | 11425.03 | 0.01 | 10524.41 | 43.02 | 2 | 11425.87 | 11425.11 | 0.01 | 10681.76 | 1.71 | 1.71 | 0.00 |
| 2-20-2 | 1 | 2 | 9733.14 | 9733.14 | 0.00 | 9044.99 | 27.11 | 2 | 9733.14 | 9733.14 | 0.00 | 9103.12 | 2.19 | 2.18 | 0.01 |
| 2-20-3 | 1 | 2 | 10096.93 | 10096.93 | 0.00 | 9590.12 | 15.39 | 2 | 10096.93 | 10096.93 | 0.00 | 9910.44 | 0.48 | 0.48 | 0.00 |
| 2-20-4 | 1 | 2 | 9495.48 | 9495.48 | 0.00 | 9136.34 | 22.02 | 2 | 9495.48 | 9495.48 | 0.00 | 9398.44 | 0.17 | 0.17 | 0.00 |
| 2-20-5 | 1 | 2 | 8302.15 | 8301.96 | 0.00 | 7805.06 | 57.90 | 2 | 8302.15 | 8301.43 | 0.01 | 7649.74 | 2.67 | 2.57 | 0.10 |
| 2-25-1 | 1 | 2 | 11440.08 | 11439.13 | 0.01 | 10733.49 | 229.65 | 2 | 11440.08 | 11439.20 | 0.01 | 11142.45 | 1.56 | 1.55 | 0.01 |
| 2-25-2 | 1 | 2 | 12429.99 | 12428.96 | 0.01 | 11747.39 | 158.96 | 2 | 12429.99 | 12429.99 | 0.00 | 12240.95 | 0.48 | 0.47 | 0.02 |
| 2-25-3 | 1 | 2 | 10977.26 | 10976.17 | 0.01 | 10216.26 | 970.95 | 2 | 10977.26 | 10976.21 | 0.01 | 10323.97 | 36.30 | 36.29 | 0.01 |
| 2-25-4 | 1 | 2 | 12279.34 | 12279.19 | 0.01 | 11705.69 | 66.51 | 2 | 12279.34 | 12279.34 | 0.00 | 11897.48 | 0.71 | 0.71 | 0.00 |
| 2-25-5 | 1 | 2 | 11791.64 | 11790.53 | 0.01 | 11101.93 | 139.22 | 2 | 11791.64 | 11791.64 | 0.00 | 11347.29 | 1.01 | 1.01 | 0.00 |
| 2-30-1 | 1 | 2 | 15001.95 | 15000.47 | 0.01 | 13939.44 | 5548.21 | 2 | 15001.95 | 15001.95 | 0.00 | 14642.40 | 9.22 | 9.20 | 0.02 |
| 2-30-2 | 1 | 2 | 12798.37 | 12797.29 | 0.01 | 11913.27 | 10161.12 | 2 | 12798.37 | 12797.30 | 0.01 | 12086.54 | 50.06 | 49.97 | 0.09 |
| 2-30-3 | 1 | 2 | 12238.97 | 12237.76 | 0.01 | 11183.54 | 7670.58 | 2 | 12238.97 | 12237.77 | 0.01 | 11754.69 | 37.15 | 36.96 | 0.20 |
| 2-30-4 | 1 | 2 | 11591.20 | 11590.04 | 0.01 | 11128.80 | 596.63 | 2 | 11591.25 | 11590.27 | 0.01 | 11320.55 | 4.53 | 4.37 | 0.16 |
| 2-30-5 | 1 | 9 | 13323.57 | 12999.39 | 2.43 | 12321.11 | 14400.00 | 2* | 13266.07* | 13264.78 | 0.01 | 12417.94 | 481.98 | 481.65 | 0.33 |
| 2-30-5 | 4 | 2 | 13266.11 | 13264.84 | 0.01 | 12321.11 | 7896.29 | 2 | 13266.11 | 13264.86 | 0.01 | 12417.94 | 103.70 | 100.70 | 3.00 |
| 2-35-1 | 4 | 9 | 14287.93 | 13909.90 | 2.65 | 13108.78 | 43200.00 | 2* | 14287.93 | 14286.59 | 0.01 | 13841.58 | 126.95 | 126.35 | 0.60 |
| 2-35-2 | 4 | 2 | 17450.55 | 17449.28 | 0.01 | 16138.00 | 3397.50 | 2 | 17449.48 | 17447.78 | 0.01 | 17072.86 | 8.11 | 8.09 | 0.02 |
| 2-35-3 | 4 | 2 | 14696.52 | 14695.11 | 0.01 | 13828.64 | 11644.95 | 2 | 14696.52 | 14696.05 | 0.00 | 14081.58 | 60.38 | 59.85 | 0.53 |
| 2-35-4 | 4 | 2 | 17694.87 | 17693.15 | 0.01 | 16599.03 | 25076.29 | 2 | 17694.87 | 17693.18 | 0.01 | 17002.00 | 105.62 | 105.45 | 0.17 |
| 2-35-5 | 4 | 2 | 16817.38 | 16815.82 | 0.01 | 15740.37 | 18041.84 | 2 | 16817.38 | 16816.08 | 0.01 | 16200.23 | 18.63 | 18.05 | 0.58 |
| 2-40-1 | 4 | 2 | 17008.84 | 17007.30 | 0.01 | 16050.14 | 21628.14 | 2 | 17008.84 | 17007.47 | 0.01 | 16568.41 | 19.32 | 19.00 | 0.32 |
| 2-40-2 | 4 | 2 | 17954.65 | 17952.87 | 0.01 | 16900.34 | 41897.40 | 2 | 17954.65 | 17953.01 | 0.01 | 17364.03 | 13.61 | 13.30 | 0.31 |
| 2-40-3 | 4 | 2 | 18085.02 | 18083.44 | 0.01 | 17039.06 | 5495.69 | 2 | 18086.39 | 18085.02 | 0.01 | 17833.45 | 5.51 | 5.32 | 0.19 |
| 2-40-4 | 4 | 9 | 18565.55 | 18164.30 | 2.16 | 17526.75 | 43200.01 | 2* | 18565.55 | 18563.70 | 0.01 | 18195.81 | 249.08 | 248.90 | 0.18 |
| 2-40-5 | 4 | 9 | 13886.66 | 13218.20 | 4.81 | 12651.37 | 43200.01 | 2* | 13804.31* | 13802.98 | 0.01 | 13115.10 | 1559.39 | 1558.00 | 1.39 |
| 2-45-1 | 4 | 2 | 18661.97 | 18660.36 | 0.01 | 17684.04 | 8093.59 | 2 | 18661.97 | 18661.97 | 0.00 | 18428.72 | 7.90 | 7.13 | 0.77 |
| 2-45-2 | 4 | 9 | 19596.73 | 19337.28 | 1.32 | 18204.88 | 43200.01 | 2* | 19596.73 | 19594.78 | 0.01 | 18476.98 | 273.86 | 270.75 | 3.11 |
| 2-45-3 | 4 | 2 | 20214.00 | 20212.55 | 0.01 | 19205.33 | 7452.02 | 2 | 20214.00 | 20212.20 | 0.01 | 19817.93 | 37.35 | 37.05 | 0.30 |
| 2-45-4 | 4 | 9 | 17322.92 | 16859.55 | 2.68 | 16154.74 | 43200.01 | 2* | 17313.46* | 17311.73 | 0.01 | 16890.72 | 672.85 | 669.75 | 3.10 |
| 2-45-5 | 4 | 9 | 24318.93 | 24075.15 | 1.00 | 22998.56 | 43200.02 | 2* | 24318.93 | 24316.59 | 0.01 | 23786.81 | 20.98 | 20.90 | 0.08 |
| 2-50-1 | 4 | 9 | 24723.83 | 24141.05 | 2.36 | 23646.84 | 43200.00 | 2* | 24706.24* | 24703.78 | 0.01 | 24153.58 | 965.72 | 961.40 | 4.32 |
| 2-50-2 | 4 | 9 | 21947.05 | 21112.52 | 3.80 | 20476.37 | 43200.03 | 2* | 21947.05 | 21944.86 | 0.01 | 21273.96 | 184.15 | 183.35 | 0.80 |
| 2-50-3 | 4 | 9 | 19707.00 | 19170.53 | 2.72 | 18708.88 | 43200.01 | 2* | 19707.00 | 19705.19 | 0.01 | 19457.56 | 95.09 | 94.05 | 1.04 |
| 2-50-4 | 4 | 9 | 19849.76 | 19707.04 | 0.72 | 19078.55 | 43200.01 | 2* | 19849.76 | 19849.63 | 0.00 | 19394.82 | 55.23 | 54.15 | 1.08 |
| 2-50-5 | 4 | 9 | 23729.09 | 23252.55 | 2.01 | 22417.39 | 43200.01 | 2* | 23729.09 | 23728.12 | 0.00 | 23333.99 | 18.49 | 18.05 | 0.44 |

**Table D.2.** The detailed results of the experiments for the objective setting R.

| Instance | # of Threads | CAR algorithm | | | | | | New formulation | | | | | | |
|---|---|---|---|---|---|---|---|---|---|---|---|---|---|---|
| | | Status | UP | LB | Bound (%) | RLB | Run Time | Status | UP | LB | Bound (%) | RLB | Run Time | Solver Time | Preprocessing Time |
| 1-10-1 | 1 | 2 | 2930.38 | 2930.38 | 0.00 | 2816.44 | 0.71 | 2 | 2930.38 | 2930.38 | 0.00 | 2754.28 | 0.09 | 0.09 | 0.00 |
| 1-10-2 | 1 | 2 | 4426.11 | 4426.11 | 0.00 | 4337.39 | 0.14 | 2 | 4426.11 | 4426.11 | 0.00 | 4426.11 | 0.06 | 0.06 | 0.00 |
| 1-10-3 | 1 | 2 | 4252.31 | 4252.31 | 0.00 | 4205.49 | 0.22 | 2 | 4252.31 | 4252.31 | 0.00 | 3871.66 | 0.07 | 0.07 | 0.00 |
| 1-10-4 | 1 | 2 | 4105.34 | 4105.34 | 0.00 | 3692.51 | 1.20 | 2 | 4105.34 | 4105.34 | 0.00 | 3600.40 | 0.19 | 0.19 | 0.00 |
| 1-10-5 | 1 | 2 | 4224.95 | 4224.95 | 0.00 | 4070.91 | 0.20 | 2 | 4224.95 | 4224.95 | 0.00 | 4164.72 | 0.07 | 0.07 | 0.00 |
| 1-15-1 | 1 | 2 | 6601.64 | 6601.64 | 0.00 | 6286.34 | 14.71 | 2 | 6601.64 | 6601.63 | 0.00 | 6420.71 | 0.49 | 0.48 | 0.02 |
| 1-15-2 | 1 | 2 | 4114.68 | 4114.68 | 0.00 | 3831.17 | 6.17 | 2 | 4114.68 | 4114.68 | 0.00 | 3804.80 | 7.11 | 6.37 | 0.74 |
| 1-15-3 | 1 | 2 | 5575.05 | 5574.72 | 0.01 | 4816.13 | 760.99 | 2 | 5575.05 | 5574.52 | 0.01 | 4973.79 | 58.04 | 57.95 | 0.09 |
| 1-15-4 | 1 | 2 | 5125.92 | 5125.42 | 0.01 | 4811.33 | 31.27 | 2 | 5125.92 | 5125.92 | 0.00 | 5004.16 | 0.67 | 0.67 | 0.00 |
| 1-15-5 | 1 | 2 | 6900.05 | 6899.40 | 0.01 | 6210.11 | 70.00 | 2 | 6900.05 | 6900.05 | 0.00 | 6311.64 | 0.96 | 0.95 | 0.01 |
| 1-20-1 | 1 | 2 | 7720.07 | 7720.07 | 0.00 | 7638.25 | 2.30 | 2 | 7720.07 | 7720.07 | 0.00 | 7652.93 | 0.28 | 0.27 | 0.01 |
| 1-20-2 | 1 | 2 | 8912.60 | 8912.25 | 0.00 | 8574.66 | 44.12 | 2 | 8912.60 | 8911.80 | 0.01 | 8738.55 | 1.72 | 1.71 | 0.01 |
| 1-20-3 | 1 | 2 | 8219.74 | 8219.74 | 0.00 | 8112.60 | 18.20 | 2 | 8219.74 | 8219.74 | 0.00 | 7993.40 | 5.73 | 5.70 | 0.03 |
| 1-20-4 | 1 | 2 | 6229.45 | 6229.02 | 0.01 | 6073.77 | 29.26 | 2 | 6229.45 | 6229.45 | 0.00 | 6155.64 | 6.90 | 6.65 | 0.25 |
| 1-20-5 | 1 | 2 | 7269.39 | 7268.70 | 0.01 | 6796.59 | 717.31 | 2 | 7269.39 | 7269.39 | 0.00 | 6925.47 | 6.70 | 6.65 | 0.05 |
| 1-25-1 | 1 | 2 | 9962.50 | 9961.54 | 0.01 | 9814.61 | 59.70 | 2 | 9962.50 | 9961.52 | 0.01 | 9494.16 | 58.06 | 57.95 | 0.11 |
| 1-25-2 | 1 | 2 | 8064.60 | 8063.79 | 0.01 | 7611.93 | 7817.43 | 2 | 8064.60 | 8064.60 | 0.00 | 7821.49 | 684.41 | 682.10 | 2.31 |
| 1-25-3 | 1 | 2 | 9422.58 | 9421.68 | 0.01 | 9182.19 | 609.38 | 2 | 9422.57 | 9422.24 | 0.00 | 9007.64 | 30.61 | 30.40 | 0.21 |
| 1-25-4 | 1 | 2 | 9460.56 | 9459.75 | 0.01 | 9135.62 | 443.76 | 2 | 9460.56 | 9460.56 | 0.00 | 9331.04 | 2.95 | 2.85 | 0.10 |
| 1-25-5 | 1 | 9 | 10398.10 | 10213.59 | 1.77 | 9786.77 | 14400.00 | 2* | 10398.10 | 10397.22 | 0.01 | 10182.71 | 26.95 | 26.60 | 0.35 |
| 1-25-5 | 4 | 2 | 10398.10 | 10397.71 | 0.00 | 9786.77 | 2049.04 | 2 | 10398.10 | 10398.10 | 0.00 | 10182.71 | 8.90 | 8.55 | 0.35 |
| 1-30-1 | 1 | 2 | 9165.77 | 9164.87 | 0.01 | 8841.70 | 2440.83 | 2 | 9165.77 | 9165.77 | 0.00 | 9027.76 | 75.22 | 61.75 | 13.47 |
| 1-30-2 | 1 | 2 | 11811.35 | 11810.92 | 0.01 | 11564.88 | 360.63 | 2 | 11811.35 | 11811.35 | 0.00 | 11745.69 | 25.28 | 23.75 | 1.53 |
| 1-30-3 | 1 | 9 | 11530.36 | 11289.64 | 2.09 | 10932.38 | 14400.00 | 2* | 11530.36 | 11529.31 | 0.01 | 11020.21 | 1179.92 | 1178.95 | 0.97 |
| 1-30-3 | 4 | 2 | 11530.36 | 11529.29 | 0.01 | 10932.38 | 28169.07 | 2 | 11530.36 | 11529.90 | 0.01 | 11020.21 | 381.92 | 380.95 | 0.97 |
| 1-30-4 | 1 | 2 | 11670.88 | 11669.72 | 0.01 | 11281.72 | 2699.82 | 2 | 11670.88 | 11669.96 | 0.01 | 11473.94 | 132.42 | 129.20 | 3.22 |
| 1-30-5 | 1 | 9 | 11330.90 | 11222.62 | 0.96 | 10924.57 | 14400.00 | 2* | 11266.00* | 11266.00 | 0.00 | 11126.19 | 28.55 | 23.75 | 4.80 |
| 1-30-5 | 4 | 2 | 11266.00 | 11264.97 | 0.01 | 10924.57 | 3958.71 | 2 | 11266.00 | 11266.00 | 0.00 | 11126.19 | 41.85 | 37.05 | 4.80 |
| 1-35-1 | 4 | 9 | 11719.16 | 11237.17 | 4.11 | 10962.85 | 43200.02 | 2* | 11583.58* | 11582.42 | 0.01 | 11113.35 | 17011.39 | 16939.45 | 71.94 |
| 1-35-2 | 4 | 9 | 12145.55 | 11927.59 | 1.80 | 11495.50 | 43200.01 | 2* | 12145.55 | 12144.54 | 0.01 | 11815.03 | 3115.54 | 3082.75 | 32.79 |
| 1-35-3 | 4 | 2 | 12346.39 | 12345.18 | 0.01 | 11955.24 | 20641.95 | 2 | 12346.39 | 12346.39 | 0.00 | 12314.08 | 175.79 | 164.35 | 11.44 |
| 1-35-4 | 4 | 2 | 12925.29 | 12924.11 | 0.01 | 12528.64 | 30214.80 | 2 | 12925.29 | 12925.29 | 0.00 | 12789.97 | 59.16 | 56.05 | 3.11 |
| 1-35-5 | 4 | 2 | 12384.35 | 12383.27 | 0.01 | 11900.52 | 29125.97 | 2 | 12384.35 | 12383.35 | 0.01 | 12240.93 | 39.86 | 36.10 | 3.76 |
| 1-40-1 | 4 | 9 | 14546.65 | 14055.27 | 3.38 | 13792.55 | 43200.01 | 2* | 14496.52* | 14495.23 | 0.01 | 14259.24 | 473.96 | 467.40 | 6.56 |
| 1-40-2 | 4 | 2 | 15741.84 | 15740.91 | 0.01 | 15181.12 | 13947.51 | 2 | 15741.84 | 15740.57 | 0.01 | 15603.98 | 86.37 | 80.75 | 5.62 |
| 1-40-3 | 4 | 9 | 13230.15 | 12735.17 | 3.74 | 12415.01 | 43200.05 | 2* | 13230.15 | 13230.15 | 0.00 | 13094.32 | 628.43 | 584.25 | 44.18 |
| 1-40-4 | 4 | 2 | 15189.51 | 15122.57 | 0.44 | 14560.88 | 43200.01 | 2* | 15189.51 | 15189.51 | 0.00 | 15122.91 | 86.93 | 76.95 | 9.98 |
| 1-40-5 | 4 | 9 | 14554.86 | 14448.54 | 0.73 | 14075.17 | 43200.01 | 2* | 14554.86 | 14553.96 | 0.01 | 14369.63 | 250.24 | 248.90 | 1.34 |
| 1-45-1 | 4 | 9 | 13832.45 | 13245.66 | 4.24 | 12999.74 | 43200.05 | 2* | 13582.63* | 13581.46 | 0.01 | 13313.76 | 33663.65 | 31279.70 | 2383.95 |
| 1-45-2 | 4 | 9 | 18423.50 | 18026.81 | 2.15 | 17692.14 | 43200.01 | 2* | 18387.82* | 18385.98 | 0.01 | 18003.43 | 2717.49 | 2665.50 | 51.99 |
| 1-45-3 | 4 | 9 | 17548.68 | 17283.63 | 1.51 | 16809.36 | 43200.01 | 2* | 17458.68 | 17546.98 | -0.51 | 17347.96 | 1497.56 | 1434.50 | 63.06 |
| 1-45-4 | 4 | 9 | 15231.54 | 14767.46 | 3.05 | 14429.80 | 43200.01 | 2* | 15188.08* | 15186.63 | 0.01 | 14938.78 | 1430.96 | 1378.45 | 52.51 |
| 1-45-5 | 4 | 9 | 17467.00 | 17126.75 | 1.95 | 16781.09 | 43200.01 | 2* | 17467.00 | 17465.59 | 0.01 | 17206.59 | 557.84 | 512.05 | 45.79 |
| 2-10-1 | 1 | 2 | 4668.58 | 4668.58 | 0.00 | 4475.21 | 0.19 | 2 | 4668.58 | 4668.58 | 0.00 | 4475.21 | 0.10 | 0.10 | 0.00 |
| 2-10-2 | 1 | 2 | 5447.17 | 5447.17 | 0.00 | 5385.90 | 0.07 | 2 | 5447.17 | 5447.17 | 0.00 | 5447.17 | 0.04 | 0.04 | 0.00 |
| 2-10-3 | 1 | 2 | 4946.66 | 4946.66 | 0.00 | 4909.98 | 0.41 | 2 | 4946.66 | 4946.66 | 0.00 | 4857.82 | 0.04 | 0.04 | 0.00 |
| 2-10-4 | 1 | 2 | 5775.16 | 5775.16 | 0.00 | 5046.24 | 0.85 | 2 | 5775.16 | 5775.16 | 0.00 | 5267.87 | 0.10 | 0.10 | 0.00 |
| 2-10-5 | 1 | 2 | 5178.47 | 5178.47 | 0.00 | 4785.74 | 0.60 | 2 | 5178.47 | 5178.47 | 0.00 | 4845.25 | 0.07 | 0.07 | 0.00 |
| 2-15-1 | 1 | 2 | 6427.57 | 6427.57 | 0.00 | 6295.31 | 1.92 | 2 | 6427.57 | 6427.57 | 0.00 | 6126.29 | 0.12 | 0.12 | 0.00 |
| 2-15-2 | 1 | 2 | 7995.39 | 7995.39 | 0.00 | 7072.49 | 10.92 | 2 | 7995.39 | 7995.39 | 0.00 | 7382.15 | 0.33 | 0.33 | 0.00 |
| 2-15-3 | 1 | 2 | 6173.88 | 6173.88 | 0.00 | 6165.34 | 0.70 | 2 | 6173.88 | 6173.88 | 0.00 | 5982.83 | 0.10 | 0.10 | 0.00 |
| 2-15-4 | 1 | 2 | 8238.88 | 8238.88 | 0.00 | 7249.54 | 5.46 | 2 | 8238.88 | 8238.88 | 0.00 | 7195.02 | 0.13 | 0.13 | 0.00 |
| 2-15-5 | 1 | 2 | 8114.11 | 8114.11 | 0.00 | 7565.10 | 6.02 | 2 | 8114.11 | 8114.11 | 0.00 | 8072.51 | 0.09 | 0.09 | 0.00 |
| 2-20-1 | 1 | 2 | 10684.58 | 10684.58 | 0.00 | 9677.25 | 46.50 | 2 | 10684.58 | 10683.95 | 0.01 | 9981.83 | 2.28 | 2.28 | 0.00 |
| 2-20-2 | 1 | 2 | 9093.23 | 9093.23 | 0.00 | 8408.78 | 19.95 | 2 | 9093.23 | 9093.23 | 0.00 | 8503.92 | 1.91 | 1.90 | 0.01 |
| 2-20-3 | 1 | 2 | 9444.03 | 9444.03 | 0.00 | 9046.59 | 14.49 | 2 | 9444.03 | 9444.03 | 0.00 | 9273.24 | 0.29 | 0.29 | 0.00 |
| 2-20-4 | 1 | 2 | 8858.50 | 8858.50 | 0.00 | 8573.31 | 12.64 | 2 | 8858.50 | 8858.50 | 0.00 | 8777.60 | 0.16 | 0.16 | 0.00 |
| 2-20-5 | 1 | 2 | 7750.71 | 7750.44 | 0.00 | 7341.92 | 41.57 | 2 | 7750.71 | 7750.71 | 0.00 | 7132.68 | 2.00 | 1.90 | 0.10 |
| 2-25-1 | 1 | 2 | 10667.20 | 10666.17 | 0.01 | 10044.21 | 182.08 | 2 | 10667.20 | 10666.37 | 0.01 | 10405.48 | 1.91 | 1.90 | 0.01 |
| 2-25-2 | 1 | 2 | 11623.07 | 11622.40 | 0.01 | 11094.80 | 78.70 | 2 | 11623.07 | 11623.07 | 0.00 | 11438.94 | 0.40 | 0.38 | 0.02 |
| 2-25-3 | 1 | 2 | 10243.91 | 10243.00 | 0.01 | 9584.05 | 609.94 | 2 | 10243.91 | 10243.10 | 0.01 | 9639.14 | 19.01 | 19.00 | 0.01 |
| 2-25-4 | 1 | 2 | 11468.38 | 11468.22 | 0.00 | 11015.48 | 44.35 | 2 | 11468.38 | 11468.30 | 0.00 | 11110.55 | 0.48 | 0.48 | 0.00 |
| 2-25-5 | 1 | 2 | 10986.41 | 10985.34 | 0.01 | 10415.81 | 96.67 | 2 | 10986.41 | 10986.41 | 0.00 | 10593.66 | 1.10 | 1.10 | 0.00 |
| 2-30-1 | 1 | 2 | 14004.58 | 14003.18 | 0.01 | 13082.92 | 6400.93 | 2 | 14004.58 | 14004.58 | 0.00 | 13670.34 | 4.77 | 4.75 | 0.02 |
| 2-30-2 | 1 | 2 | 11958.16 | 11956.98 | 0.01 | 11162.09 | 5706.98 | 2 | 11958.16 | 11957.71 | 0.01 | 11278.54 | 48.54 | 48.45 | 0.09 |
| 2-30-3 | 1 | 2 | 11432.84 | 11431.77 | 0.01 | 10551.97 | 6124.27 | 2 | 11432.83 | 11432.83 | 0.00 | 10975.03 | 23.95 | 23.75 | 0.20 |
| 2-30-4 | 1 | 2 | 10827.64 | 10826.62 | 0.01 | 10491.07 | 372.99 | 2 | 10827.64 | 10827.64 | 0.00 | 10576.78 | 3.96 | 3.80 | 0.16 |
| 2-30-5 | 1 | 2 | 12382.16 | 12380.93 | 0.01 | 11560.96 | 13526.95 | 2 | 12382.16 | 12380.97 | 0.01 | 11578.86 | 502.88 | 502.55 | 0.33 |
| 2-35-1 | 4 | 9 | 13338.11 | 13000.93 | 2.53 | 12340.16 | 43200.02 | 2* | 13338.11 | 13337.43 | 0.01 | 12920.35 | 171.60 | 171.00 | 0.60 |
| 2-35-2 | 4 | 2 | 16285.18 | 16285.71 | 0.00 | 15114.16 | 1755.40 | 2 | 16285.18 | 16284.26 | 0.01 | 15939.98 | 17.12 | 17.10 | 0.02 |
| 2-35-3 | 4 | 2 | 13709.31 | 13707.94 | 0.01 | 13091.27 | 8732.74 | 2 | 13709.30 | 13708.11 | 0.01 | 13136.50 | 71.78 | 71.25 | 0.53 |
| 2-35-4 | 4 | 2 | 16511.45 | 16510.08 | 0.01 | 15590.82 | 9765.14 | 2 | 16511.45 | 16509.94 | 0.01 | 15859.19 | 32.47 | 32.30 | 0.17 |
| 2-35-5 | 4 | 2 | 15706.05 | 15704.77 | 0.01 | 14796.76 | 9491.46 | 2 | 15706.05 | 15706.05 | 0.00 | 15123.76 | 11.98 | 11.40 | 0.58 |
| 2-40-1 | 4 | 2 | 15866.07 | 15865.05 | 0.01 | 15050.78 | 32162.29 | 2 | 15866.07 | 15864.64 | 0.01 | 15457.95 | 26.92 | 26.60 | 0.32 |
| 2-40-2 | 4 | 9 | 16751.86 | 16563.74 | 1.12 | 15858.93 | 43200.01 | 2* | 16751.86 | 16751.7 | 0.00 | 16195.19 | 18.36 | 18.05 | 0.31 |
| 2-40-3 | 4 | 2 | 16867.91 | 16866.92 | 0.01 | 16156.49 | 3298.58 | 2 | 16867.91 | 16867.67 | 0.00 | 16645.03 | 3.04 | 2.85 | 0.19 |
| 2-40-4 | 4 | 9 | 17327.69 | 16959.05 | 2.13 | 16506.21 | 43200.02 | 2* | 17327.69 | 17326.1 | 0.01 | 16981.23 | 307.98 | 307.80 | 0.18 |
| 2-40-5 | 4 | 9 | 12989.88 | 12333.83 | 5.05 | 11903.81 | 43200.01 | 2* | 12890.2* | 12888.9 | 0.01 | 12239.55 | 1649.64 | 1648.25 | 1.39 |
| 2-45-1 | 4 | 2 | 17395.52 | 17393.95 | 0.01 | 16659.27 | 6142.20 | 2 | 17395.52 | 17395.52 | 0.00 | 17208.05 | 10.27 | 9.50 | 0.77 |
| 2-45-2 | 4 | 2 | 18299.79 | 18298.04 | 0.01 | 17418.54 | 41018.03 | 2 | 18299.78 | 18298.29 | 0.01 | 17238.43 | 235.86 | 232.75 | 3.11 |
| 2-45-3 | 4 | 9 | 18854.89 | 18627.89 | 1.20 | 18056.72 | 43200.01 | 2* | 18854.89 | 18853.70 | 0.01 | 18473.82 | 79.15 | 78.85 | 0.30 |
| 2-45-4 | 4 | 9 | 16156.67 | 15830.19 | 2.02 | 15187.75 | 43200.01 | 2* | 16156.67 | 16155.06 | 0.01 | 15763.46 | 724.15 | 721.05 | 3.10 |
| 2-45-5 | 4 | 2 | 22695.80 | 22693.71 | 0.01 | 21613.50 | 37956.46 | 2 | 22695.80 | 22695.74 | 0.00 | 22200.83 | 139.73 | 139.65 | 0.08 |
| 2-50-1 | 4 | 9 | 23045.51 | 22630.72 | 1.80 | 22117.36 | 43200.01 | 2* | 23045.51 | 23043.23 | 0.01 | 22533.05 | 979.97 | 975.65 | 4.32 |
| 2-50-2 | 4 | 9 | 20734.05 | 19945.00 | 3.81 | 19248.06 | 43200.01 | 2* | 20472.32* | 20470.29 | 0.01 | 19853.91 | 389.35 | 388.55 | 0.80 |
| 2-50-3 | 4 | 9 | 18482.61 | 17942.60 | 2.92 | 17595.49 | 43200.01 | 2* | 18399.27* | 18399.19 | 0.00 | 18164.82 | 93.19 | 92.15 | 1.04 |
| 2-50-4 | 4 | 9 | 18537.57 | 18293.14 | 1.32 | 17778.54 | 43200.01 | 2* | 18537.57 | 18535.74 | 0.01 | 18104.46 | 70.43 | 69.35 | 1.08 |
| 2-50-5 | 4 | 9 | 22160.84 | 21761.23 | 1.80 | 20922.11 | 43200.02 | 2* | 22137.56* | 22135.64 | 0.01 | 21781.83 | 25.14 | 24.70 | 0.44 |

**Table D.3.** The detailed results of the experiments for the objective setting E.

| Instance | # of Threads | CAR algorithm | | | | | | New formulation | | | | | | |
|---|---|---|---|---|---|---|---|---|---|---|---|---|---|---|
| | | Status | UP | LB | Bound (%) | RLB | Run Time | Status | UP | LB | Bound (%) | RLB | Run Time | Solver Time | Preprocessing Time |
| 1-10-1 | 1 | 2 | 203.26 | 203.26 | 0.00 | 187.01 | 0.75 | 2 | 203.26 | 203.26 | 0.00 | 190.93 | 0.15 | 0.15 | 0.00 |
| 1-10-2 | 1 | 2 | 314.61 | 314.61 | 0.00 | 313.09 | 0.44 | 2 | 314.61 | 314.61 | 0.00 | 311.92 | 0.07 | 0.07 | 0.00 |
| 1-10-3 | 1 | 2 | 305.50 | 305.50 | 0.00 | 303.89 | 0.48 | 2 | 305.50 | 305.50 | 0.00 | 278.98 | 0.12 | 0.12 | 0.00 |
| 1-10-4 | 1 | 2 | 287.37 | 287.37 | 0.00 | 252.62 | 1.86 | 2 | 287.37 | 287.37 | 0.00 | 258.60 | 0.31 | 0.31 | 0.00 |
| 1-10-5 | 1 | 2 | 301.05 | 301.05 | 0.00 | 288.76 | 1.07 | 2 | 301.05 | 301.05 | 0.00 | 297.70 | 0.09 | 0.09 | 0.00 |
| 1-15-1 | 1 | 2 | 473.11 | 473.09 | 0.01 | 452.17 | 17.95 | 2 | 473.11 | 473.11 | 0.00 | 461.09 | 0.57 | 0.55 | 0.02 |
| 1-15-2 | 1 | 2 | 284.54 | 284.54 | 0.00 | 254.98 | 18.54 | 2 | 284.54 | 284.54 | 0.00 | 266.98 | 4.71 | 3.97 | 0.74 |
| 1-15-3 | 1 | 2 | 387.75 | 387.73 | 0.01 | 291.80 | 318.25 | 2 | 387.75 | 387.75 | 0.00 | 360.76 | 6.74 | 6.65 | 0.09 |
| 1-15-4 | 1 | 2 | 363.64 | 363.62 | 0.01 | 319.23 | 22.92 | 2 | 363.64 | 363.64 | 0.00 | 354.14 | 0.58 | 0.58 | 0.00 |
| 1-15-5 | 1 | 2 | 484.30 | 484.30 | 0.00 | 402.73 | 35.56 | 2 | 484.30 | 484.30 | 0.00 | 445.53 | 1.29 | 1.28 | 0.01 |
| 1-20-1 | 1 | 2 | 566.03 | 566.03 | 0.00 | 552.75 | 6.19 | 2 | 566.03 | 566.03 | 0.00 | 550.09 | 0.29 | 0.29 | 0.01 |
| 1-20-2 | 1 | 2 | 638.28 | 638.22 | 0.01 | 554.84 | 35.40 | 2 | 638.28 | 638.28 | 0.00 | 627.93 | 0.90 | 0.89 | 0.01 |
| 1-20-3 | 1 | 2 | 596.95 | 596.90 | 0.01 | 525.35 | 93.63 | 2 | 596.95 | 596.89 | 0.01 | 578.41 | 3.59 | 3.56 | 0.03 |
| 1-20-4 | 1 | 2 | 459.17 | 459.17 | 0.00 | 379.52 | 141.39 | 2 | 459.17 | 459.17 | 0.00 | 431.49 | 4.77 | 4.52 | 0.25 |
| 1-20-5 | 1 | 2 | 505.48 | 505.44 | 0.01 | 481.93 | 84.89 | 2 | 505.48 | 505.48 | 0.00 | 488.57 | 1.11 | 1.06 | 0.05 |
| 1-25-1 | 1 | 2 | 712.27 | 712.20 | 0.01 | 636.29 | 76.87 | 2 | 712.27 | 712.25 | 0.00 | 689.29 | 4.22 | 4.10 | 0.11 |
| 1-25-2 | 1 | 2 | 575.15 | 575.09 | 0.01 | 496.12 | 10468.16 | 2 | 575.15 | 575.15 | 0.00 | 558.56 | 194.11 | 191.81 | 2.31 |
| 1-25-3 | 1 | 2 | 675.22 | 675.16 | 0.01 | 562.38 | 1244.17 | 2 | 675.22 | 675.17 | 0.01 | 650.94 | 16.80 | 16.59 | 0.21 |
| 1-25-4 | 1 | 2 | 689.22 | 689.18 | 0.01 | 581.51 | 736.09 | 2 | 689.22 | 689.22 | 0.00 | 673.71 | 3.83 | 3.72 | 0.10 |
| 1-25-5 | 1 | 2 | 747.82 | 747.75 | 0.01 | 614.92 | 10017.91 | 2 | 747.82 | 747.75 | 0.01 | 726.14 | 23.44 | 23.09 | 0.35 |
| 1-30-1 | 1 | 2 | 667.91 | 667.85 | 0.01 | 554.16 | 3224.96 | 2 | 667.91 | 667.91 | 0.00 | 650.78 | 98.02 | 84.55 | 13.47 |
| 1-30-2 | 1 | 2 | 857.98 | 857.89 | 0.01 | 734.00 | 942.35 | 2 | 857.98 | 857.98 | 0.00 | 853.38 | 17.17 | 15.64 | 1.53 |
| 1-30-3 | 1 | 9 | 831.79 | 807.84 | 2.88 | 706.68 | 14400.01 | 2* | 827.98* | 827.98 | 0.00 | 797.28 | 181.47 | 180.50 | 0.97 |
| 1-30-3 | 4 | 2 | 827.98 | 827.90 | 0.01 | 706.68 | 9770.65 | 2 | 827.98 | 827.98 | 0.00 | 797.28 | 41.63 | 40.66 | 0.97 |
| 1-30-4 | 1 | 2 | 845.53 | 845.44 | 0.01 | 755.72 | 11981.76 | 2 | 845.52 | 845.52 | 0.00 | 832.59 | 12.50 | 9.28 | 3.22 |
| 1-30-5 | 1 | 9 | 814.68 | 798.34 | 2.01 | 729.81 | 14400.01 | 2* | 814.68 | 814.68 | 0.00 | 798.83 | 40.14 | 35.34 | 4.80 |
| 1-30-5 | 4 | 9 | 814.68 | 800.26 | 1.77 | 729.81 | 43200.01 | 2* | 814.68 | 814.68 | 0.00 | 798.83 | 39.00 | 34.20 | 4.80 |
| 1-35-1 | 4 | 9 | 846.23 | 807.70 | 4.55 | 691.89 | 43200.01 | 2* | 844.27* | 844.19 | 0.01 | 805.41 | 4476.14 | 4404.20 | 71.94 |
| 1-35-2 | 4 | 9 | 879.95 | 848.50 | 3.58 | 735.60 | 43200.01 | 2 | 879.96 | 879.96 | 0.00 | 859.71 | 419.44 | 386.65 | 32.79 |
| 1-35-3 | 4 | 2 | 889.46 | 889.38 | 0.01 | 762.91 | 21866.89 | 2 | 889.46 | 889.46 | 0.00 | 884.65 | 190.04 | 178.60 | 11.44 |
| 1-35-4 | 4 | 2 | 941.71 | 941.62 | 0.01 | 808.55 | 36895.74 | 2 | 941.71 | 941.71 | 0.00 | 920.76 | 59.16 | 56.05 | 3.11 |
| 1-35-5 | 4 | 2 | 894.50 | 894.41 | 0.01 | 752.61 | 15121.73 | 2 | 894.50 | 894.50 | 0.00 | 883.41 | 37.96 | 34.20 | 3.76 |
| 1-40-1 | 4 | 9 | 1050.10 | 1019.41 | 2.92 | 899.93 | 43200.01 | 2* | 1047.06* | 1047.00 | 0.01 | 1023.79 | 78.76 | 72.20 | 6.56 |
| 1-40-2 | 4 | 9 | 1144.12 | 1133.85 | 0.90 | 1036.83 | 43200.00 | 2* | 1144.12 | 1144.00 | 0.01 | 1132.13 | 144.32 | 138.70 | 5.62 |
| 1-40-3 | 4 | 9 | 960.94 | 915.43 | 4.74 | 759.84 | 43200.02 | 2* | 947.76* | 947.76 | 0.00 | 940.41 | 567.63 | 523.45 | 44.18 |
| 1-40-4 | 4 | 2 | 1099.03 | 1098.92 | 0.01 | 984.30 | 23208.12 | 2 | 1099.03 | 1098.92 | 0.01 | 1088.53 | 93.58 | 83.60 | 9.98 |
| 1-40-5 | 4 | 9 | 1056.13 | 1042.42 | 1.30 | 996.06 | 43200.00 | 2* | 1056.12 | 1056.03 | 0.01 | 1040.60 | 172.34 | 171.00 | 1.34 |
| 1-45-1 | 4 | 9 | 991.73 | 955.58 | 3.65 | 840.48 | 43200.01 | 2* | 984.19* | 984.10 | 0.01 | 957.23 | 35632.05 | 33248.10 | 2383.95 |
| 1-45-2 | 4 | 9 | 1343.02 | 1307.24 | 2.66 | 1203.82 | 43200.01 | 2* | 1334.75* | 1334.68 | 0.01 | 1312.94 | 7993.99 | 7942.00 | 51.99 |
| 1-45-3 | 4 | 9 | 1284.48 | 1236.21 | 3.76 | 1053.93 | 43200.01 | 2* | 1278.69* | 1278.58 | 0.01 | 1256.83 | 2508.36 | 2445.30 | 63.06 |
| 1-45-4 | 4 | 9 | 1104.80 | 1055.93 | 4.42 | 928.92 | 43200.01 | 2* | 1103.09 | 1102.98 | 0.01 | 1074.92 | 923.66 | 871.15 | 52.51 |
| 1-45-5 | 4 | 9 | 1254.39 | 1231.99 | 1.79 | 1116.21 | 43200.01 | 2* | 1254.42 | 1254.31 | 0.01 | 1233.81 | 843.79 | 798.00 | 45.79 |
| 2-10-1 | 1 | 2 | 324.19 | 324.19 | 0.00 | 292.51 | 0.46 | 2 | 324.19 | 324.19 | 0.00 | 307.32 | 0.21 | 0.21 | 0.00 |
| 2-10-2 | 1 | 2 | 379.88 | 379.88 | 0.00 | 378.69 | 0.13 | 2 | 379.88 | 379.88 | 0.00 | 377.32 | 0.07 | 0.07 | 0.00 |
| 2-10-3 | 1 | 2 | 324.87 | 324.87 | 0.00 | 323.63 | 0.24 | 2 | 324.87 | 324.87 | 0.00 | 316.70 | 0.06 | 0.06 | 0.00 |
| 2-10-4 | 1 | 2 | 382.90 | 382.90 | 0.00 | 311.92 | 1.34 | 2 | 382.89 | 382.89 | 0.00 | 350.05 | 0.18 | 0.18 | 0.00 |
| 2-10-5 | 1 | 2 | 356.87 | 356.87 | 0.00 | 309.05 | 0.36 | 2 | 356.87 | 356.87 | 0.00 | 328.96 | 0.11 | 0.11 | 0.00 |
| 2-15-1 | 1 | 2 | 444.24 | 444.24 | 0.00 | 428.56 | 1.15 | 2 | 444.24 | 444.24 | 0.00 | 419.47 | 0.14 | 0.14 | 0.00 |
| 2-15-2 | 1 | 2 | 533.99 | 533.99 | 0.00 | 446.37 | 7.60 | 2 | 533.99 | 533.99 | 0.00 | 501.44 | 0.49 | 0.49 | 0.00 |
| 2-15-3 | 1 | 2 | 436.30 | 436.30 | 0.00 | 398.59 | 3.11 | 2 | 436.30 | 436.30 | 0.00 | 418.88 | 0.11 | 0.11 | 0.00 |
| 2-15-4 | 1 | 2 | 541.30 | 541.30 | 0.00 | 463.68 | 6.52 | 2 | 541.30 | 541.30 | 0.00 | 482.04 | 0.14 | 0.14 | 0.00 |
| 2-15-5 | 1 | 2 | 543.98 | 543.98 | 0.00 | 508.07 | 2.07 | 2 | 543.98 | 543.98 | 0.00 | 538.87 | 0.11 | 0.11 | 0.00 |
| 2-20-1 | 1 | 2 | 741.26 | 741.26 | 0.01 | 643.61 | 57.14 | 2 | 741.29 | 741.23 | 0.01 | 699.93 | 1.62 | 1.62 | 0.00 |
| 2-20-2 | 1 | 2 | 631.21 | 631.21 | 0.00 | 540.84 | 36.43 | 2 | 631.21 | 631.20 | 0.00 | 597.78 | 1.53 | 1.52 | 0.01 |
| 2-20-3 | 1 | 2 | 652.89 | 652.89 | 0.00 | 577.35 | 21.07 | 2 | 652.89 | 652.89 | 0.00 | 636.30 | 0.35 | 0.35 | 0.00 |
| 2-20-4 | 1 | 2 | 619.46 | 619.46 | 0.00 | 595.00 | 4.57 | 2 | 619.46 | 619.45 | 0.00 | 614.06 | 0.20 | 0.20 | 0.00 |
| 2-20-5 | 1 | 2 | 551.44 | 551.44 | 0.00 | 463.92 | 67.24 | 2 | 551.44 | 551.40 | 0.01 | 511.55 | 3.52 | 3.42 | 0.10 |
| 2-25-1 | 1 | 2 | 752.17 | 752.10 | 0.01 | 646.18 | 570.93 | 2 | 752.17 | 752.17 | 0.00 | 730.89 | 0.67 | 0.67 | 0.01 |
| 2-25-2 | 1 | 2 | 805.33 | 805.25 | 0.01 | 707.52 | 186.40 | 2 | 805.33 | 805.33 | 0.00 | 788.23 | 0.49 | 0.48 | 0.02 |
| 2-25-3 | 1 | 2 | 717.81 | 717.74 | 0.01 | 601.74 | 668.59 | 2 | 717.81 | 717.81 | 0.00 | 679.78 | 12.74 | 12.73 | 0.01 |
| 2-25-4 | 1 | 2 | 805.73 | 805.73 | 0.01 | 778.82 | 24.42 | 2 | 805.73 | 805.73 | 0.00 | 786.64 | 0.34 | 0.34 | 0.00 |
| 2-25-5 | 1 | 2 | 784.06 | 783.99 | 0.01 | 668.29 | 54.81 | 2 | 784.06 | 784.06 | 0.00 | 752.18 | 0.57 | 0.57 | 0.00 |
| 2-30-1 | 1 | 2 | 983.41 | 983.31 | 0.01 | 851.98 | 1125.56 | 2 | 983.41 | 983.37 | 0.00 | 964.88 | 7.52 | 7.51 | 0.02 |
| 2-30-2 | 1 | 9 | 840.21 | 814.08 | 3.11 | 701.91 | 14400.00 | 2* | 840.21 | 840.21 | 0.00 | 807.32 | 9.59 | 9.50 | 0.09 |
| 2-30-2 | 4 | 2 | 840.21 | 840.15 | 0.01 | 701.91 | 1641.33 | 2 | 840.21 | 840.19 | 0.00 | 807.32 | 9.32 | 7.79 | 1.53 |
| 2-30-3 | 1 | 2 | 802.56 | 802.48 | 0.01 | 660.53 | 12469.98 | 2 | 802.56 | 802.49 | 0.01 | 777.24 | 45.80 | 45.60 | 0.20 |
| 2-30-4 | 1 | 2 | 763.56 | 763.48 | 0.01 | 645.44 | 7057.37 | 2 | 763.56 | 763.56 | 0.00 | 741.41 | 3.87 | 3.71 | 0.16 |
| 1-30-5 | 1 | 9 | 878.70 | 855.74 | 2.61 | 739.77 | 14400.00 | 2* | 878.70 | 878.62 | 0.01 | 827.88 | 99.42 | 99.09 | 0.33 |
| 1-30-5 | 4 | 2 | 878.70 | 878.64 | 0.01 | 739.77 | 30463.13 | 2 | 878.70 | 878.70 | 0.00 | 827.88 | 38.15 | 35.15 | 3.00 |
| 2-35-1 | 4 | 2 | 941.88 | 941.78 | 0.01 | 804.14 | 31873.37 | 2 | 941.88 | 941.81 | 0.01 | 911.61 | 48.77 | 48.17 | 0.60 |
| 2-35-2 | 4 | 2 | 1158.52 | 1158.43 | 0.01 | 986.13 | 18306.00 | 2 | 1158.52 | 1158.41 | 0.01 | 1122.98 | 11.61 | 11.59 | 0.02 |
| 2-35-3 | 4 | 9 | 982.43 | 954.59 | 2.83 | 822.41 | 43200.01 | 2* | 982.43 | 982.37 | 0.01 | 940.28 | 43.09 | 42.56 | 0.53 |
| 2-35-4 | 4 | 2 | 1171.34 | 1171.25 | 0.01 | 997.27 | 38648.37 | 2 | 1171.34 | 1171.23 | 0.01 | 1137.53 | 30.38 | 30.21 | 0.17 |
| 2-35-5 | 4 | 2 | 1108.78 | 1108.68 | 0.01 | 946.20 | 42292.75 | 2 | 1108.78 | 1108.78 | 0.00 | 1073.35 | 9.03 | 8.46 | 0.58 |
| 2-40-1 | 4 | 2 | 1131.77 | 1131.67 | 0.01 | 982.10 | 6219.51 | 2 | 1131.77 | 1131.77 | 0.00 | 1106.27 | 7.82 | 7.51 | 0.32 |
| 2-40-2 | 4 | 9 | 1189.64 | 1172.57 | 1.44 | 1010.32 | 43200.01 | 2* | 1189.64 | 1189.64 | 0.00 | 1157.67 | 9.33 | 9.03 | 0.31 |
| 2-40-3 | 4 | 2 | 1202.43 | 1202.31 | 0.01 | 1023.70 | 30554.31 | 2 | 1202.42 | 1202.40 | 0.00 | 1186.04 | 4.69 | 4.50 | 0.19 |
| 2-40-4 | 4 | 2 | 1233.52 | 1233.41 | 0.01 | 1111.12 | 7308.37 | 2 | 1233.52 | 1233.41 | 0.01 | 1203.90 | 121.78 | 121.60 | 0.18 |
| 2-40-5 | 4 | 9 | 920.01 | 861.24 | 6.39 | 763.25 | 43200.01 | 2* | 900.96* | 900.88 | 0.01 | 860.28 | 131.54 | 130.15 | 1.39 |
| 2-45-1 | 4 | 2 | 1242.93 | 1242.93 | 0.00 | 1076.34 | 1802.46 | 2 | 1242.93 | 1242.93 | 0.00 | 1215.02 | 9.29 | 8.52 | 0.77 |
| 2-45-2 | 4 | 9 | 1288.56 | 1271.78 | 1.30 | 1085.49 | 43200.02 | 2* | 1288.65 | 1288.56 | 0.01 | 1232.97 | 126.61 | 123.50 | 3.11 |
| 2-45-3 | 4 | 9 | 1355.49 | 1326.35 | 2.15 | 1169.61 | 43200.01 | 2* | 1355.49 | 1355.36 | 0.01 | 1327.11 | 28.80 | 28.50 | 0.30 |
| 2-45-4 | 4 | 9 | 1142.63 | 1121.53 | 1.85 | 979.63 | 43200.01 | 2* | 1142.61 | 1142.61 | 0.00 | 1120.95 | 50.60 | 47.50 | 3.10 |
| 2-45-5 | 4 | 9 | 1617.45 | 1583.33 | 2.11 | 1403.58 | 43200.01 | 2* | 1617.45 | 1617.30 | 0.01 | 1581.75 | 15.28 | 15.20 | 0.08 |
| 2-50-1 | 4 | 9 | 1649.17 | 1637.05 | 0.74 | 1513.33 | 43200.01 | 2* | 1649.17 | 1649.09 | 0.01 | 1618.29 | 67.97 | 63.65 | 4.32 |
| 2-50-2 | 4 | 9 | 1461.93 | 1429.13 | 2.24 | 1271.46 | 43200.01 | 2* | 1457.63* | 1457.52 | 0.01 | 1410.30 | 75.85 | 75.05 | 0.80 |
| 2-50-3 | 4 | 9 | 1304.48 | 1249.69 | 4.20 | 1149.86 | 43200.01 | 2* | 1304.48 | 1304.48 | 0.00 | 1286.16 | 63.74 | 62.70 | 1.04 |
| 2-50-4 | 4 | 9 | 1309.09 | 1266.09 | 3.29 | 1128.84 | 43200.10 | 2* | 1306.60* | 1306.60 | 0.00 | 1284.93 | 70.43 | 69.35 | 1.08 |
| 2-50-5 | 4 | 9 | 1578.18 | 1553.65 | 1.55 | 1454.56 | 43200.01 | 2* | 1578.18 | 1578.14 | 0.00 | 1540.25 | 13.74 | 13.30 | 0.44 |